\documentclass[11pt]{extarticle}
\usepackage[english]{babel}
\usepackage{graphicx}
\usepackage{framed}
\usepackage[normalem]{ulem}
\usepackage{amsmath}
\usepackage{amsthm}
\usepackage{amssymb}
\usepackage{amsfonts}
\usepackage{mathtools}
\usepackage{commath}
\usepackage{enumerate}
\usepackage[utf8]{inputenc}
\usepackage{bbm}
\usepackage{xcolor}
\usepackage[top=1 in,bottom=1in, left=1 in, right=1 in]{geometry}
\usepackage{enumitem}
\usepackage[symbol]{footmisc}
\usepackage{url}
\usepackage[linesnumbered,boxruled]{algorithm2e}

\usepackage{lineno}

\theoremstyle{plain}
\newtheorem{theorem}{Theorem}
\newtheorem{proposition}{Proposition}
\newtheorem{lemma}{Lemma}
\newtheorem{corollary}{Corollary}
\newtheorem{assumption}{Assumption}
\newtheorem*{definition}{Definition}
\newtheorem{rem}{Remark}

\newcommand{\R}{\mathbb{R}}

\newcommand{\M}{\mathcal{M}}

\usepackage{color}
\definecolor{darkgreen}{rgb}{0.0, 0.2, 0.13}
\newcommand{\rot}[1]{{\color{black} #1}}

\newcommand{\vth}{\vartheta}

\newcommand{\taum}{\tau_\text{m}}
\newcommand{\taus}{\tau_\text{s}}

\newcommand{\lrrund}[1]{\left( #1 \right)}
\newcommand{\lreckig}[1]{\left[ #1 \right]}

\title{On a finite-size neuronal population equation}
\author{Valentin Schmutz\thanks{Brain Mind Institute, École Polytechnique Fédérale de Lausanne, 1015 Lausanne, Switzerland} \and Eva Löcherbach\thanks{Statistique, Analyse et Modélisation Multidisciplinaire,  EA 4543 et FR FP2M 2036 CNRS, Université Paris 1 Panthéon-Sorbonne, 75013 Paris, France} \and Tilo Schwalger\thanks{Institut für Mathematik, Technische Universität Berlin, 10623 Berlin, Germany}\,\,\,\thanks{Bernstein Center for Computational Neuroscience Berlin, 10115 Berlin, Germany}}
\date{}

\begin{document}

\maketitle

\begin{abstract}
Population equations for infinitely large networks of spiking neurons have a long tradition in theoretical neuroscience. In this work, we analyze a recent generalization of these equations to populations of finite size, which takes the form of a nonlinear stochastic integral equation.  We prove that,  in the case of leaky integrate-and-fire (LIF) neurons with escape noise and for a  slightly simplified version of the model, the equation is well-posed and stable in the sense of Brémaud-Massoulié. 
The proof combines methods from Markov processes taking values in the space of positive measures and nonlinear Hawkes processes. For applications, we also provide efficient simulation algorithms.
\end{abstract}

\vspace{0.5cm}


\section{Introduction}
Neuronal population equations describe the dynamics of large networks of neurons in terms of single neuron parameters \cite{GerKis14}. As such, they are useful mathematical abstractions for relating microscopic and large-scale brain signals, and contribute to the biophysical interpretation of the latter \cite{DecJir08}. Their motivation is twofold: on the one hand, they enable the theoretical analysis of emergent phenomena, like collective oscillations \cite{BruHak99, Ger00, CorTan21}; on the other hand, from the data analysis point of view, they constitute the basis of `forward models' of large-scale brain signals \cite{DecJir08, SanKno15,CaiIye16, Bre17, FriPre19}. This second motivation requires neuronal population equations to achieve the right balance between accuracy (the equation faithfully captures the dynamics of the population of neurons it represents) and usability (the equation can be efficiently simulated).

An example of such neuronal population equation is the integral equation (or refractory density equation) for a homogeneous network of spiking neurons (``neuronal population'') \cite{Ger95,Ger00,ChiGra07,GerKis14, SchChi19}. Contrary to standard neural-mass models \cite{WilCow72,DecJir08,JanRit95}, the integral equation captures the effect of neuronal refractoriness on the mean population dynamics \cite{ChiGra07,GerKis14,SchChi19}, and is exact in the mean-field limit if neurons are modeled as intensity-based renewal point processes \cite{DemGal15,FouLoe16,Che17}. Specific examples of the integral equation are the time-elapsed neuron network model \cite{PakPer10} (or age-structured model \cite{DumHen16b}) and the voltage-structured model of \cite{DemGal15,FouLoe16}.

Besides capturing the effect of single neuron dynamics (such as post-spike refractory effects) on the \emph{mean} population dynamics, there is a second challenge for neuronal population equations: the proper account of \emph{fluctuations}. Fluctuations of the average population activity arise in the case of finite population sizes, and vanish in the mean-field limit of infinitely many neurons. From a modeling perspective, an important question arises: Are the relevant neuronal populations large enough so that finite-size fluctuations can be neglected? There is no clear answer to this question but the anatomical and functional organization of the cerebral cortex into different cortical areas, columns and layers each containing different cell classes \cite{HarShe15,PotDie14,SchBak18,BilCai20} requires a subdivision of a cortical circuit into many, relatively small populations. For example, at the scale of a cortical column, empirical data from mouse barrel cortex suggests populations of around 50 to 2000 neurons \cite{LefTom09}. For these  population sizes, finite-size fluctuations are non-negligible and this noise may strongly impact the nonlinear population dynamics \cite{SchDeg17}. Therefore, modeling cortical circuits at the mesoscopic scale of populations requires a stochastic description, in marked contrast to the deterministic integral equation.

Rigorous extensions of the integral equations to account for finite-size fluctuations are subject to an accuracy/usability trade-off. If neuronal refractoriness is neglected, the population equation reduces to that of \cite{DelFou16,DitLoe17} and finite-size noise can be added, by the linear-noise approximation \cite{HeeSta21}, or granting some Markov embedding, by the diffusion approximation \cite{DitLoe17}, whose numerical implementation is relatively simple \cite{CheMel20}. These approaches fail to reproduce the non-stationary dynamics of the mean population activity and the temporal correlation structure of fluctuations for a population of spiking neurons with refractoriness (Fig.~\ref{fig:1}a). On the other hand, if one does not neglect refractoriness, central limit theorem-based arguments lead to formal SPDE's \cite{Che17b,DumPay17}, which are computationally expensive to simulate, or to formal integral equations with colored noise \cite{DegSch14}, for which a simulation algorithm is unknown. 

In \cite{SchDeg17}, a heuristic extension of the integral equation with finite-size fluctuations is derived. It can be easily simulated and takes into account the effects of neuronal refractoriness. While this extension is not exact, its numerical implementation gives an accurate approximation to the dynamics of finite-size networks of spiking neurons, such as the broad class of generalized integrate-and-fire neurons \cite{PozMen15,SchDeg17} and formal renewal-type neurons \cite{Ger00,PieGal20}. Moreover, since it takes the form of an intensity-based point process, the likelihood of a population spike train can be easily computed, which enables efficient data fitting \rot{\cite{RenLon20,WanSch22}}. The intensity function of this point process exhibits a novel type of nonlinear history dependence that goes beyond nonlinear Hawkes processes and has not been studied mathematically so far. In particular, the stability of the process observed in simulations is poorly understood from the theoretical point of view.
Therefore, the aim of this work is to give a rigorous foundation to the model of \cite{SchDeg17} and prove its stability. 

Below, we briefly give a review of some standard population equations. We then present the finite-size model of \cite{SchDeg17} in a slightly simplified form. Finally, we show that the simplified model, in the case of leaky integrate-and-fire (LIF) neurons with escape noise \cite{Ger00,GalLoe16}, can be written as a SPDE driven by Poisson noise, which will be the main object of study in this work.

\subsection{Neuronal population equations}\label{sec:review}
To give a mathematical introduction to the integral equation formalism, it is useful to consider the special case of LIF neurons with escape noise \cite{Ger00,GalLoe16}, which is also the main case we will treat in this work. Let us consider a network of $N$ identical neurons that are all-to-all connected with uniform connection strength $J/N$ for $J\in\R$. Each neuron $i$ has a voltage variable $U^{i,N}$ which evolves according to the system of SDE's:
For all $i=1, \dots, N$,
\begin{linenomath}\begin{subequations}\label{eq:network_GL}
  \begin{align}
    dU_t^{i,N}&=\frac{\mu_t-U_t^{i,N}}{\taum}dt-U_{t^-}^{i,N}dZ_t^{i,N}+\frac{J}{N}\sum_{j=1}^NdZ_t^{j,N},\label{eq:network_GL_U}\\ 
    Z_t^{i,N}&=\int_{[0,t]\times\R_+} \mathbbm{1}_{z\leq f(U_{s^-}^{i,N})}\pi^i(ds,dz).
  \end{align}
\end{subequations}\end{linenomath}
Here,  $Z_t^{i,N}$ is the spike counting process of the neuron $i$ and has intensity $f(U_{t^-}^{i,N})$, $t^-$ denoting the left limit. Furthermore, $\mu_t$ comprises the resting potential and the (possibly time-dependent) external drive, $\taum$ is the membrane time constant, $f:\R \to\R_+$ is the intensity function and $\{\pi^i\}_{i=1,\dots,N}$ is a collection of independent Poisson random measures on $\R_+\times\R_+$ with Lebesgue intensity measure.

Equation \eqref{eq:network_GL} is called a \emph{microscopic model} because the neuronal dynamics is modeled with single-cell resolution (Fig.~\ref{fig:1}a, top).  A drastic reduction of the complexity of the model can be achieved by coarse-graining over the population of neurons. To this end, we consider the \emph{empirical population activity}
\begin{linenomath}\begin{equation}
  \label{eq:empir-A}
  A_{t,\mathfrak{h}}^N=\frac{1}{N}\sum_{i=1}^N\frac{Z_{t+\mathfrak{h}}^{i,N}-Z_t^{i,N}}{\mathfrak{h}},
\end{equation}\end{linenomath}
where $\mathfrak{h}>0$ is a small time interval determining the temporal resolution (Fig.~\ref{fig:1}a, bottom). Neuronal population equations are models of such coarse-grained quantitities that describe the neuronal dynamics at the scale of whole populations. If the population is of finite size ($N<\infty$), the dynamics is called a \emph{mesoscopic model}, while the dynamics for an infinitely large population ($N\rightarrow\infty$) is referred to as a \emph{macroscopic model}.
 In \cite{DemGal15,FouLoe16}, the authors proved that in the macroscopic limit $N\to\infty$, if the initial conditions $\{U_0^i\}_{i=1,\dots,N}$ are $i.i.d.$ with law $\nu_0$, the empirical measure of the system~\eqref{eq:network_GL} is characterized  by the voltage-structured PDE (with solutions in the sense of measures \cite{CorTan21}): 
For all $u\in\R$ and $t>0$,
\begin{subequations}\label{eq:voltage_PDE}
\begin{linenomath}\begin{align}
    &\partial_t \rho(du,t) + \partial_u\left(\left(\frac{\mu_t-u}{\taum} + J\rho_t[f]\right)\rho(du,t)\right) = -f(u)\rho(du,t)+\rho_t[f]\delta_0(du), \label{eq:voltage_PDE_continuity}\\
    &\rho_0 = \nu_0,
\end{align}\end{linenomath}
\end{subequations}
where $\rho_t := \rho(\cdot, t)$ and $\rho_t[f] := \int_{\R} f(u)\rho(du,t)$.

The latter can be interpreted as the population activity
\begin{equation}
  \label{eq:pop-act-Ninfty}
  \lim_{\mathfrak{h}\downarrow 0}\lim_{N\rightarrow\infty}A_{t,\mathfrak{h}}^N=A(t):=\rho_t[f].
\end{equation}
Furthermore, $\rho_t[1]=1$ for all $t>0$ expressing the fact that the number of neurons is conserved.

We now transform Eq.~\eqref{eq:voltage_PDE} into an integral equation. For all continuous functions $a:\R_+\to\R$, we define the time-dependent vector field $b^{a}(t,u):= (\mu_t-u)/\taum + J{a}(t)$ and write, for all $0\leq s\leq t$, $\Phi^{a}_{s,t}(u)$ the associated flow given by 
\begin{equation}\label{eq:Phia}
    \Phi^{a}_{s,t}(u):= u e^{-\frac{t-s}{\taum}} + \int_s^te^{-\frac{t-r}{\taum}}\frac{\mu_r}{\taum}\,dr + J\int_s^t e^{-\frac{t-r}{\taum}}a_r d r , \qquad \forall u\in\R .
\end{equation}
We can now define, for all $0\leq s\leq t$, 
\begin{equation}\label{eq:lambda_A}
    \lambda^a(t|s):=f(\Phi^a_{s,t}(0)) \qquad \text{and} \qquad S^a(t|s):= \exp\left(-\int_s^t \lambda^a(r|s)\,dr\right).
\end{equation}
The function $\lambda^a(t|s)$, called hazard rate, gives the intensity at time $t$ (i.e. the instantaneous probability of emitting a spike) as a function of the time of the last spike $s$ and the past population activity $(a(r))_{s\leq r\leq t}$; the membrane potential dynamics of LIF neurons -- leaky integration and spike-triggered reset, Eq.~\eqref{eq:network_GL_U} -- are accounted for in the definition of $\lambda^a(t|s)$. Similarly, the function $S^a(t|s)$, called the survival, gives the probability of not emitting a spike in the time interval $]s,t[$, given that the last spike was emitted at time $s$. By the method of characteristics, we get that the population activity $A(t)$ solves the integral equation
\begin{equation}\label{eq:integral}
    A(t) = H^A(t) + \int_0^t \lambda^A(t|s)S^A(t|s)A(s)ds, 
\end{equation}
where 
\begin{equation}\label{eq:H}
    {H^A(t):= \int_{\R} f(\Phi^A_{0,t}(u))e^{-\int_0^t f(\Phi^A_{0,r}(u))dr}\nu_0(du).}
\end{equation} 
Equation~\eqref{eq:integral} is the integral equation of \cite{WilCow72,Ger95,Ger00}, see also \cite{CorTan20}. Note that, traditionally, the integral equation has no explicit initial condition and therefore requires a normalizing condition \cite[Sec.~14.1]{GerKis14}. The integral equation~\eqref{eq:integral} is normalized such that
\begin{equation}
  \label{eq:normalization}
  \widetilde{H}^A(t)+\int_0^t S^A(t|s)A(s)\,ds=1  
\end{equation}
for all $t>0$, where we defined
\begin{equation}
  \label{eq:H_tilde}
  {\widetilde{H}^{A}(t):= \int_{\R} e^{-\int_0^t f(\Phi^{A}_{0,r}(u))dr}\nu_0(du).}
\end{equation}
The normalization, Eq.~\eqref{eq:normalization}, expresses the fact that the number of neurons is conserved.\footnote{The conservation of neuronal mass can be understood as follows: At time $t$, $\widetilde{H}^A(t)$ represents the fraction of neurons (\#neurons divided by $N$) that had their unique last spike before time $0$, while for $s\in[0,t[$ the term $S^A(t|s)A(s)ds$ represents the fraction of neurons that had their unique last spike time in the interval $[s,s+ds[$ (here $A(s)ds$ is the fraction of neurons that fired in that interval and $S^A(t|s)$ is the probability for one neuron of not emetting a spike in $]s,t[$ given a spike at time $s$). Therefore, $\int_0^tS^A(t|s)A(s)\,ds$ represents the fraction of neurons that had their unique last spike in $[0,t[$. Hence, Eq.~\eqref{eq:normalization} states that the fraction of neurons at time $t$ that had their unique last spike time before time $t$ (either before time $0$ or since time $0$) is equal to unity. Since this statement holds for all $t>0$ and each neuron has exactly one last spike time before time $t$, the total number of neurons must be conserved.} \rot{Note that the integral equation~\eqref{eq:integral} is simply the time derivative of the normalizing condition Eq.~\eqref{eq:normalization}; this fact has been originally used to derive the integral equation \cite{Ger00}.}

In the case of LIF neurons with escape noise, the voltage-structured equation~\eqref{eq:voltage_PDE} is equivalent to the integral equation~\eqref{eq:integral} if $\lambda^A(t|s)$ is defined by Eq.~\eqref{eq:lambda_A}. However, we could have chosen a different definition for the hazard rate $\lambda^A(t|s)$; the integral equation is therefore more general than Eq.~\eqref{eq:voltage_PDE}. In fact, Eq.~\eqref{eq:integral} can be seen as a renewal equation that holds for any population of neurons modeled as time-inhomogeneous renewal processes \cite{PieGal20}. For example, the Fokker-Planck equation for neuronal networks with diffusive noise (see \cite[Ch.~13]{GerKis14}) or the time-elapsed neuron network model \cite{PakPer10} can also be written as an integral equation with a suitable choice of the hazard rate.

\subsection{The finite-size integral equation}\label{sec:finite_size_integral_equation}

\begin{figure}[p]
  \centering
  \includegraphics[width=\textwidth]{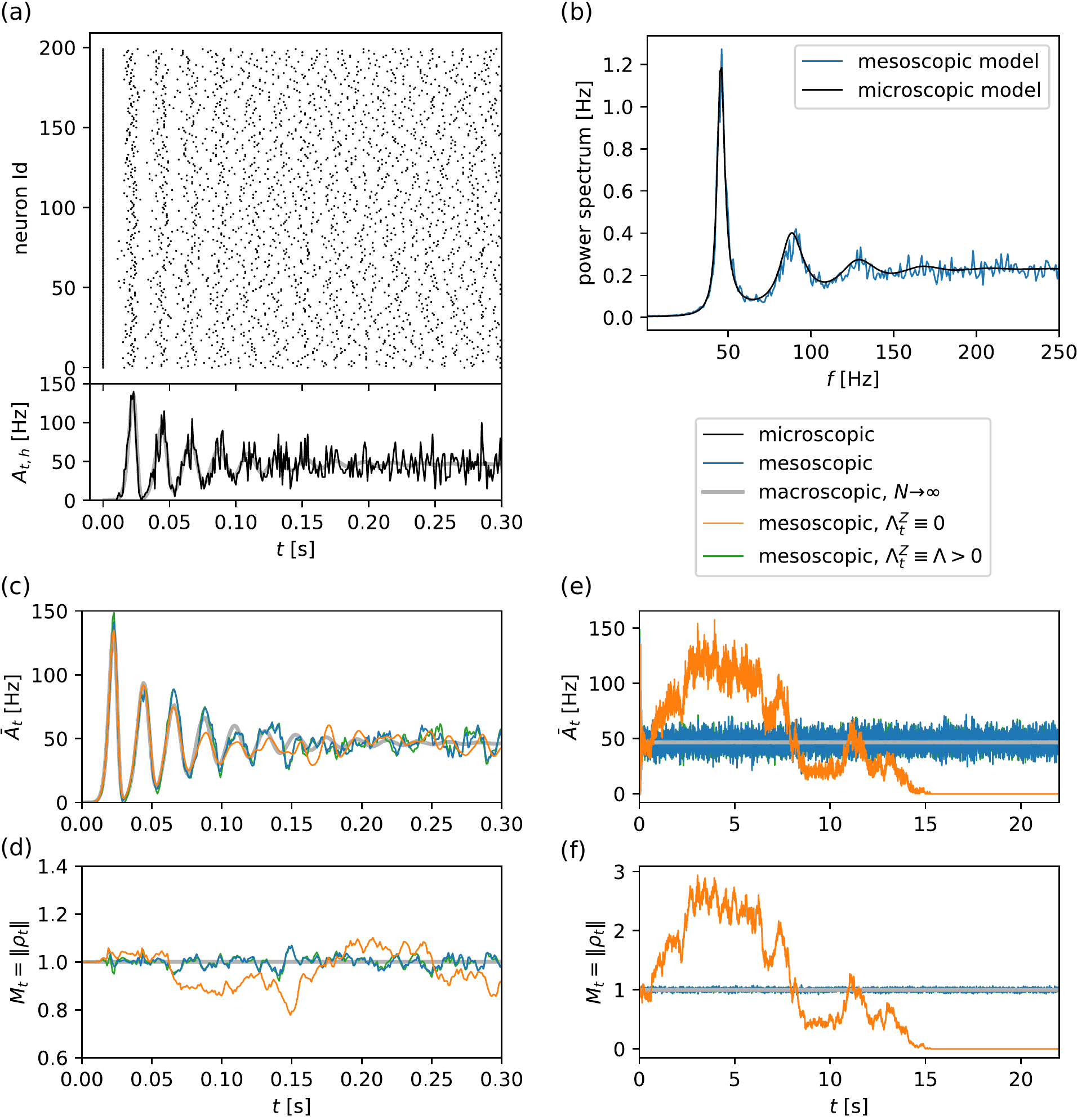}
  \caption{Mesoscopic population dynamics. (a) Top: Spike-raster plot of a microscopic model of $N=200$ uncoupled LIF neurons with escape noise, Eq.~\eqref{eq:network_GL} with $J=0$. Neurons were initialized in a synchronized state, i.e. all neurons spiked at time $t=0$. Bottom: Empirical population activity measured with temporal bin size $\mathfrak{h}=0.001$s (black line) and macroscopic population activity predicted by the deterministic integral equation \eqref{eq:integral} for $N\to\infty$ with $\nu_0=\delta_0$ (gray line). (b) Comparison of the power spectral densities (as defined in Sec.~\ref{sec:psd}, see also \cite{SchDeg17}) 
      of the empirical population activities $A_{t,\mathfrak{h}}(t)$ of the microscopic model (black line, exact theory \cite{GerKis14}) and $\hat{A}_{t,\mathfrak{h}}(t)$ of the mesoscopic model (blue line, simulation). (c,d) $\bar{A}_t$ \eqref{eq:stochastic_integra_A_bar} and mass $M_t$ \eqref{eq:mass} for simulations of the mesoscopic model (blue line) and the `naive' mesoscopic model with $\Lambda^Z_t \equiv 0$ (orange line). For comparison, the macroscopic model and the mesoscopic model with fixed $\Lambda_t^Z\equiv 277$~Hz (corresponding to the temporally averaged $\Lambda_t^Z$ of the mesoscopic model) are shown by gray and green lines, respectively. (e,f) Same as (c,d) but for a longer simulation time. Parameters: $\taum=0.02$~s, $\mu=20$~mV, $f(u)=ce^{(u-\vth)/\Delta_u}$, $c=10$~Hz, $\vth=10$~mV, $\Delta_u=1$~mV.}
  \label{fig:1}
\end{figure}

In \cite{SchDeg17}, the authors derive a generalization of the integral equation~\eqref{eq:integral} which takes into account finite-size noise. For clarity, we will present the equation of \cite{SchDeg17} in the case of LIF neurons with escape noise. Before presenting the model, we need to extend the definitions Eq.~\eqref{eq:lambda_A}. For all non-decreasing functions $z:\R_+ \ni t\mapsto z_t$ with bounded variation on finite time intervals, we redefine, for all $0\leq s\leq t$,
\begin{equation}\label{eq:Phiz}
    \Phi^{z}_{s,t}(u):= u e^{-\frac{t-s}{\taum}} + \int_s^te^{-\frac{t-r}{\taum}}\frac{\mu_r}{\taum}\,dr + J\int_{]s, t ]}  e^{-\frac{t-r}{\taum}}dz_r, \qquad \forall u\in\R .
\end{equation}
We can now extend the definitions Eqs.~\eqref{eq:lambda_A}, \eqref{eq:H} and \eqref{eq:H_tilde}, replacing $\Phi^{A}$ by Eq.~\eqref{eq:Phiz}.

For a finite number of neurons $N$, the finite-size integral equation of \cite{SchDeg17} (``mesoscopic model'') can be written as follows: For all $t\geq 0$,
\begin{subequations}\label{eq:stochastic_integral}
\begin{linenomath}\begin{align}
    Z_t &= \frac{1}{N}\int_{[0,t]\times \R_+}\mathbbm{1}_{z\leq N\bar{A}_{s^-}}\pi(ds,dz),\\
    \bar{A}_t &= \left[H^{Z}(t)+\int_{[0, t ]} \lambda^{Z}(t|s) S^{Z}(t|s) dZ_s+ \Lambda^Z_t\left(1-\widetilde{H}^{Z}(t)-\int_{[0,t]}  S^{Z}(t|s)dZ_s\right)\right]_+, \label{eq:stochastic_integra_A_bar} \\
    \Lambda^Z_t &= \frac{G^Z(t) + \int_{[0,t]} \lambda^Z(t|s)\{1 - S^Z(t|s)\}S^Z(t|s)dZ_s}{\widetilde{G}^Z(t) + \int_{[0,t]} \{1 - S^Z(t|s)\}S^Z(t|s)dZ_s},\label{eq:stochastic_integra_Lambda}
\end{align}\end{linenomath}
\end{subequations}
where $\pi$ is a Poisson random measure on $\R_+\times\R_+$ with Lebesgue intensity measure and $[\,\cdot\,]_+ = \max(0,\cdot)$. The functions $G^Z$ and $\widetilde{G}^Z$ are analogous to $H^Z$ and $\widetilde{H}^Z$:
\begin{linenomath}\begin{align*}
    G^Z(t) &:= \int_{\R} f(\Phi^Z_{0,t}(u))\left\{1 - e^{-\int_0^t f(\Phi^Z_{0,r}(u))dr}\right\}e^{-\int_0^t f(\Phi^Z_{0,r}(u))dr}\nu_0(du),\\
    \widetilde{G}^{Z}(t) &:= \int_{\R} \left\{1 - e^{-\int_0^t f(\Phi^Z_{0,r}(u))dr}\right\}e^{-\int_0^t f(\Phi^Z_{0,r}(u))dr}\nu_0(du).
\end{align*}\end{linenomath}
The mesoscopic model, Eq.~\eqref{eq:stochastic_integral} defines a jump process $Z_t$ where jumps of size $1/N$ occur with intensity $N\bar{A}_{t^-}$. The derivation of Eq.~\eqref{eq:stochastic_integral}, explained in detail in \cite[pp.~35--43]{SchDeg17}, involves heuristic arguments and approximations. Consequently, this mesoscopic model is inexact (in contrast to the formal SPDE of \cite{Che17b, DumPay17}). However, extensive numerical simulations have shown that the model is highly accurate in many multiscale modeling applications \cite{SchDeg17} (see also Figure~\ref{fig:1}b). Moreover, it has the advantage of being an intensity-based and history-dependent point process, and as such, can be efficiently simulated and used for statistical data analysis \cite{RenLon20}. A concise derivation of Eq.~\eqref{eq:stochastic_integral} is presented in Section~\ref{sec:link}, where we also show that for some convenient initial condition, the functions $H^{Z}$, $\widetilde{H}^{Z}$, $G^{Z}$ and $\widetilde{G}^{Z}$ are trivial.

The finite-size analog of the population activity $A(t)$ for infinitely large populations (Eq.~\eqref{eq:integral}) is the distributional derivative of $Z_t$,
  \begin{linenomath}\begin{equation*}
    \dot{Z}_t= \frac{1}{N}\sum_{k}\delta(t-{t_{k}}),
  \end{equation*}\end{linenomath}
where $t_k$ are the jump times of $Z_t$ and $\delta(\cdot)$ denotes the Dirac delta distribution\footnote[2]{Formally, $\dot{Z}_t dt:=dZ_t$, where $dZ$ is the Lebesgue-Stieltjes measure associated with the counting measure $Z$.}. We call $\dot Z_t$  the population spike train (sum of $\delta$-pulses at spike times $t_{k}$). Note that the biologically relevant quantity is the empirical population activity at a finite time resolution, $\hat{A}_{t,\mathfrak{h}}^N:=\mathfrak{h}^{-1}\int_t^{t+\mathfrak{h}^+}\dot Z_s\,ds=[Z(t+\mathfrak{h})-Z(t)]/\mathfrak{h}$, for some small time interval $\mathfrak{h}>0$.   Furthermore, we will often call the finite-size population model, Eq.~\eqref{eq:stochastic_integral}, \emph{mesoscopic model} in contrast to ``macroscopic model'' that refers to the case $N\rightarrow\infty$. Note that the variables $\bar{A}_t$ and $Z_t$ describe the neuronal activity of the population as a whole, driven by only one single Poisson noise $\pi(dt,[0,N\bar{A}_{t^-}])$. A time discretization of the mesoscopic model permits an efficient simulation of the neuronal dynamics directly on the population level, without the need to simulate individual neurons (see Sec.~\ref{sec:simu} and Algorithm~\ref{algo}). Importantly, even though the mesoscopic model is an approximation, it accurately captures the  statistics of the population activity $A_{t,h}^N$ of the original microscopic model. In particular, the fluctuation statistics of the population activities $A_{t,h}^N$ and $\hat{A}_{t,h}^N$, as expressed by their power spectral density, are well matched (Fig.~\ref{fig:1}b, also see \cite{SchDeg17} for further examples). 


A key difference between the macroscopic model for an infinitely large population, Eq.~\eqref{eq:integral},  and the mesoscopic model, Eq.~\eqref{eq:stochastic_integral}, is the `correction term' $\Lambda^Z_t(1-\dots)$ in Eq.~\eqref{eq:stochastic_integra_A_bar} arising due to finite network size, $N<\infty$. This correction term may seem unexpected in light of the following heuristic argument: in Eq.~\eqref{eq:integral} for infinite $N$, the fraction of neurons $A(s)ds$  firing in the past, $s<t$, contribute to the current activity $A(t)dt$ with probability $\lambda^A(t|s)S^A(t|s)dt$. For finite $N$, the corresponding fraction of neurons is $dZ_s$, and assuming that the probability to fire their next spike at time $t$ is again given by $\lambda^A(t|s)S^A(t|s)dt$, the expected activity should be given by the much simpler expression $\bar A_{t,\text{naive}}=H^{Z}(t)+\int_0^t \lambda^{Z}(t|s) S^{Z}(t|s) dZ_s$. This naive finite-size model is obtained by putting $\Lambda^Z_t\equiv 0$, and thus lacks the 'correction term'. Numerical simulations of the naive finite-size model indeed reproduce the transient initial dynamics of the population activity at short times, including damped oscillations caused by refractoriness (Fig.~\ref{fig:1}c, orange curve). However, longer simulations of the naive model reveal that the population rate $\bar A_t$ strongly fluctuates and eventually collapses to the silent solution $\bar A_t=0$. In contrast, the mesoscopic model, Eq.~\eqref{eq:stochastic_integral} with $\Lambda_t^Z>0$,  reaches a non-silent, stationary state consistent with  the microscopic model, Eq.~\eqref{eq:network_GL} (Fig.~\ref{fig:1}e). A completely open theoretical question is: Why does the `correction term' in Eq.~\eqref{eq:stochastic_integra_A_bar} `stabilize' the finite-size neuronal population dynamics? 

To address this question mathematically, we focus our analysis on the case where the modulating factor $\Lambda^Z_t$ is fixed ($\Lambda^Z_t\equiv\Lambda>0$). This is a simplified version of the finite-size integral equation~\eqref{eq:stochastic_integral}, for which we can prove a rigorous stability result. Note that fixing $\Lambda^Z_t\equiv\Lambda>0$ is for mathematical tractability only; for practical modeling, $\Lambda^Z_t$ as defined in Eq.~\eqref{eq:stochastic_integra_Lambda} should be preferred (a detailed simulation algorithm is presented in Section~\ref{sec:simu}). 

Before presenting our main stability result in Sec.~\ref{sec:main}, we provide some additional insights into the mechanisms of the finite-size integral equation~\eqref{eq:stochastic_integral}, in particular, why is the naive model ($\Lambda=0$) expected to fail. First, in Sec.~\ref{sec:relat-with-nonl}, we show a close relationship between the finite-size integral equation and nonlinear Hawkes processes, for which stability properties are well known. Second, in Sec.~\ref{sec:neur-mass}, we propose a heuristic argument for the stability in terms of neuronal mass conservation and an analogy with the Cox-Ingersoll-Ross process.

\subsubsection{Relationship with nonlinear Hawkes processes}
\label{sec:relat-with-nonl}

If $J=0$ (neurons do not interact), $\mu_t\equiv \mu$ (the external drive is constant) and $\Lambda^Z_t\equiv\Lambda$, Eq.~\eqref{eq:stochastic_integral} reduces to a nonlinear Hawkes process \cite{BreMas96}: For all $t\geq 0$,
\begin{subequations}\label{eq:nonlinear_hawkes}
\begin{linenomath}\begin{align}
    Z_t &= \frac{1}{N}\int_{[0,t]\times\R_+}\mathbbm{1}_{z\leq N\bar{A}_{s^-}}\pi(ds,dz),\\
    \bar{A}_t &= \Bigg[\Lambda + H^0(t) - \Lambda\widetilde{H}^0(t)+\int_{[0,t]} \underbrace{(\lambda^0(t|s)-\Lambda)S^{0}(t|s)}_{=:h^{\Lambda}(t-s)} dZ_s\Bigg]_+, \label{eq:nonlinear_hawkes_A}
\end{align}\end{linenomath}
\end{subequations}
where $\lambda^0$, $S^0$, $H^0$ and $\tilde{H}^0$ correspond to the definitions Eqs.~\eqref{eq:lambda_A}, \eqref{eq:H} and \eqref{eq:H_tilde} when $\Phi^Z$ (Eq.~\eqref{eq:Phiz}) is replaced by $\Phi^0_{s,t}(u) = u e^{-\frac{t-s}{\taum}} + \int_s^te^{-\frac{t-r}{\taum}}\frac{\mu_r}{\taum}\,dr$.

The function $h^{\Lambda}:\R_+\to\R$ in Eq.~\eqref{eq:nonlinear_hawkes_A} can be interpreted as the self-interaction kernel of the nonlinear Hawkes process. The model \eqref{eq:nonlinear_hawkes} is not particularly useful in practice since it only approximates the dynamics of a population of non-interacting neurons with constant external input. Nevertheless it sheds light on the role of $\Lambda$ on the stability of the mesoscopic model and it helps to see why the theory of nonlinear Hawkes processes \cite{BreMas96} will prove to be instrumental in this work. It is easy to verify that $\int_0^\infty h^{\Lambda}(t)dt=1$ if $\Lambda =0$ and $\int_0^\infty h^{\Lambda}(t)dt<1$ if $\Lambda > 0$. If $\Lambda = 0$, Eq.~\eqref{eq:nonlinear_hawkes} is a critical Hawkes process and has a nontrivial stationary solution only if $h^0$ is heavy-tailed \cite{BreMas01} (which is not the case for the neuron models considered here). On the other hand, if $\Lambda > 0$, Eq.~\eqref{eq:nonlinear_hawkes} is a stable nonlinear Hawkes process with a unique stationary solution (Theorem~1 in \cite{BreMas96} and see also \cite{CosGra20}). Hence, in the time-homogeneous ($\mu_t\equiv\mu$) and non-interacting case ($J=0$), $\Lambda_t\equiv\Lambda>0$ is a sufficient condition for the stability of Eq.~\eqref{eq:nonlinear_hawkes}, in the sense of \cite{BreMas96}.

To generalize this stability result to the interacting case ($J\neq0$), we will use a Markov embedding of Eq.~\eqref{eq:stochastic_integral} and Meyn-Tweedie theory \cite{MeyTwe93}, in addition to standard techniques for nonlinear Hawkes processes \cite{BreMas96}.

\subsubsection{Approximate conservation of neuronal mass}
\label{sec:neur-mass}

In contrast to the conservation of neuronal mass in the macroscopic model, Eq.~\eqref{eq:normalization}, such a strict conservation law does no longer hold for the mesoscopic model, Eq.~\eqref{eq:stochastic_integral}. However, in analogy to Eq.~\eqref{eq:normalization}, we would expect the neuronal ``mass''
\begin{equation}\label{eq:mass}
    M_t:=\widetilde{H}^Z(t)+\int_{[0, t]} S^Z(t|s) dZ_s
\end{equation}
to stay close to $1$. This feature is supported by simulations of the mesoscopic model showing that $M_t$ fluctuates around unity (Fig.~\ref{fig:1}d,f). Indeed, the number of neurons in the system~\eqref{eq:network_GL} being obviously constant, the finite-size population model Eq.~\eqref{eq:stochastic_integral} should reflect this mass conservation principle. 

Let us consider the first hitting time $\tau^*=\inf\{t>0:\bar A_t=0\}$. For $0<t<\tau^*$, the intensity $\bar{A}_t$ is strictly positive, hence Eq.~\eqref{eq:stochastic_integra_A_bar} can always be written as
\begin{linenomath}\begin{equation*}
    \bar{A}_t = H^Z(t)+\int_{[0, t]} \lambda^{Z}(t|s) S^{Z}(t|s) dZ_s+ \Lambda^Z_t\left(1-M_t\right).
\end{equation*}\end{linenomath}
By formal differentiation of Eq.~\eqref{eq:mass}, we obtain for $0<t<\tau^*$
\begin{linenomath}\begin{align}
dM_t=-H^Z(t)dt+dZ_t-\lrrund{\int_{[0, t]} \lambda^Z(t|s)S^Z(t|s)\,dZ_s}dt&=\Lambda^Z_t(1-M_t)dt+d\tilde Z_t, \label{eq:OUP}
\end{align}\end{linenomath}
where $\tilde Z_t:=Z_t - \int_0^t\bar{A}_s\,ds$ is the compensated jump process. Equation \eqref{eq:OUP} yields some rough insights into the dynamics of the neuronal mass $M_t$. For simplicity, let us assume $\Lambda_t^Z\equiv\Lambda$ to be constant. First, the conditional mean $\bar M_t^c:=\mathbb{E}[M_t|\tau^*>t]$ can be obtained by averaging Eq.~\eqref{eq:OUP}: $d\bar M_t^c=\Lambda(1-\bar M_t^c)dt$. This equation shows that its solution, $\bar M_t^c=1+(\widetilde H^Z(0)-1)e^{-\Lambda t}$, is attracted to unity if $\Lambda>0$. Conversely, in the naive model, when $\Lambda=0$, the conditional mean does not drift towards unity but remains constant, $\bar M_t^c=\widetilde H^Z(0)$ for all $t>0$. 
Second, in the naive model ($\Lambda=0$), once $M_t$ hits the boundary $0$, it sticks to this boundary forever, i.e. $M_t=0$ for all $t>\tau^*$ (Fig.~\ref{fig:1}f). In fact, if $f$ is upper bounded by $\|f\|_\infty<\infty$, we have $0\le \bar A_t\le\|f\|_\infty M_t+\Lambda(1-M_t)$. Thus, $M_t=0$ and $\Lambda=0$ entails that $\bar A_t=0$, and hence the ``noise''  $d\tilde{Z}$ in Eq.~\eqref{eq:OUP} vanishes. 

Third, if the jumps of $\tilde Z_t$ are small and frequent enough and if the increments of $\tilde Z_t$ are `independent' enough, we may replace $d\tilde{Z}_t$ by its diffusion approximation $\sqrt{\bar A_t/N}dW_t$, where $W_t$ is a Wiener process. If we further assume that $\bar A_t$ and $M_t$ vary roughly in proportion (as suggested by Fig.~\ref{fig:1}e,f for the naive model), we expect that $M_t$ behaves like a Cox-Ingersoll-Ross process, $d\hat{M}_t=\Lambda(1-\hat{M}_t)dt+\sigma\sqrt{\hat{M}_t}dW_t$ where $\sigma$ is the volatility parameter. Due to the drift term, this process fluctuates around its mean $\mathbb{E}[\hat{M}_t]=1$ if $\Lambda>0$, consistent with simulations of the model (Fig.~\ref{fig:1}d,f). Such drift force is absent in the naive model, $\Lambda=0,$ in which case $d\hat{M}_t = \sigma\sqrt{\hat{M}_t}dW_t $ describes the critical Feller branching diffusion which goes extinct in the long run (and once it hits $0$ remains there forever), with extinction probability $ P (\hat{M}_t  = 0 | \hat{M}_0 = x ) = e^{ - \frac{x}{\sigma^2 t } }.$

\subsection{Markov embedding of the finite-size integral equation}

As the voltage-structured equation~\eqref{eq:voltage_PDE} can be transformed into an integral equation, assuming $\Lambda^Z_t\equiv \Lambda$, we can transform the stochastic integral equation~\eqref{eq:stochastic_integral} back into a voltage-structured SPDE driven by Poisson noise. Denoting $\mathcal{M}_+$ the space of nonnegative finite measures on $\R$, for all $\mathcal{M}_+$-valued random variables $\widehat{\nu_0}$, the SPDE formally writes: 

For all $t>0$ and $u\in\R$,
\begin{subequations}\label{eq:SPDE}
\begin{linenomath}\begin{align}
    &\partial_t \rho(du,t) +\partial_u\left(\left(\frac{\mu_t-u}{\taum} + J \dot{Z}_t\right)\rho(du,t^-)\right) = -f(u)\rho(du,t)+\dot Z_t\delta_0(du), \\
    &Z_t = \frac{1}{N}\int_{[0,t]\times\R_+}\mathbbm{1}_{z\leq N\bar{A}_{s^-}}\pi(ds,dz) \quad \text{with }\bar{A}_t := [\rho_t[f] + \Lambda(1-\norm{\rho_t})]_+,\label{eq:def_Z}\\
    &\rho_0 = \widehat{\nu_0},
\end{align}\end{linenomath}
\end{subequations}
where $\norm{\cdot}$ denotes the total variation norm, that is, the total mass of the measure. 

We will give a precise meaning to the SPDE~\eqref{eq:SPDE} and show that it is equivalent to the stochastic integral equation \eqref{eq:stochastic_integral} in Section~\ref{sec:well-posedness} below. The two jump terms $ \partial_u (J \dot{Z}_t \rho(du,t^-))$ and $\dot Z_t\delta_0(du) $ have the following interpretation. At each jump time of $ Z_t,$ the current mass of the solution $ \rho (du, t)$ is shifted by $ J/N$ and a mass $ (1/N) \delta_0$ is added to the current value of the solution (emulating the membrane potential reset of LIF neurons, Eq.~\eqref{eq:network_GL_U}). Although the jump intensity $ N \bar A_{t^-}$ of $ Z_t$ is not a priori bounded, we shall prove in Lemma \ref{lemma:foster_lyapunov} below that almost surely $Z$ has only a finite number of jumps within each finite time interval such that Eq.~\eqref{eq:SPDE} is well-posed as a measure-valued piecewise deterministic Markov process having c\`adl\`ag trajectories. 

We say that Eq.~\eqref{eq:SPDE} is the Markov embedding of the jump process Eq.~\eqref{eq:stochastic_integral} (with $\Lambda_t^Z \equiv \Lambda$) and  that $Z$ is the jump process associated with the solution $\rho$.

\subsection{Assumptions and main result}\label{sec:main}

The main result of this work concerns the stability of Eq.~\eqref{eq:SPDE}. We use a notion of stability that is close to that of Br\'emaud and Massouli\'e \cite{BreMas96} for nonlinear Hawkes processes.

We say that a jump process $Z$ is \textit{stationary} if, for all $\tau>0$, the time-shifted process $(Z_{t+\tau}-Z_{\tau})_{t\geq 0}$ has the same law as $(Z_t - Z_0)_{t\geq 0}$. Then, we say that a solution $\bar{\rho}$ to Eq.~\eqref{eq:SPDE} with the $\M_+$-valued random initial condition $\bar{\nu}_0$ is stationary if the associated jump process $\bar{Z}$ is stationary. 

Since the noise in Eq.~\eqref{eq:SPDE} comes from a Poisson random measure, we can naturally construct a \textit{coupling} of two solutions $\rho$ and $\tilde{\rho}$ to Eq.~\eqref{eq:SPDE} (for different, possibly random, initial conditions) on the same probability space, using the same underlying Poisson random measure. Writing $Z$ and $\tilde{Z}$ the jump processes associated with $\rho$ and $\tilde{\rho}$, we define $T_c$ the \textit{coupling time} of $Z$ and $\tilde{Z}$, i.e. the time starting from which $Z$ and $\tilde{Z}$ are identical:
\begin{equation}\label{eq:tc}
    T_c := \inf\left\{\tau\geq 0 : (Z_{t+\tau}-Z_{\tau})_{t\geq 0} \equiv (\tilde{Z}_{t+\tau}-\tilde{Z}_{\tau})_{t\geq 0}\right\},
\end{equation}
with the usual convention that $T_c = +\infty$ if $Z$ and $\tilde{Z}$ never couple. In other words, $T_c$ is the time starting from which $\rho$ and $\tilde{\rho}$ have the exact same jump times. By abuse of terminology, we will say that $T_c$ is the coupling time of $\rho$ and $\tilde{\rho}$ although it is in fact the coupling time of the associated jump processes. We can now adapt the definition of stability in variation of \cite{BreMas96}:
\begin{definition}[Stability in variation]
The voltage-structured SPDE~\eqref{eq:SPDE} is \emph{stable in variation} if there exists a stationary process $\{\bar{\rho}, \bar{\nu}_0\}$ solving Eq.~\eqref{eq:SPDE} such that for all $\M_+$-valued random initial conditions $\widehat{\nu_0}$, there exists a coupling of $\bar{\rho}$ and $\rho$ (a solution to Eq.~\eqref{eq:SPDE} with initial condition $\widehat{\nu_0}$), such that the coupling time $T_c$ of $\bar{\rho}$ and $\rho$ is almost surely finite.
\end{definition}

In modeling terms, the stability in variation implies that for any (random) initial condition $\widehat{\nu_0}$, the population spike train $\dot{Z}_t$ relaxes to a unique stationary process in finite time. More specifically, for any initial condition $\widehat\nu_0\in\mathcal{M}_+$, if we draw $\bar\nu_0$ from a stationary distribution and if we simulate the two corresponding processes with the same Poisson noise, they couple in finite time almost surely. In particular, this implies the uniqueness of the stationary distribution.

To prove that Eq.~\eqref{eq:SPDE} is stable in variation, we need
\begin{assumption}\label{assumption:constant}
$\mu_t \equiv \mu \in \R$.
\end{assumption}
This just means that the external drive is time-homogeneous and it is a natural assumption to make if we want to show relaxation to a stationary process.

The other important assumption concerns the intensity function $f$:
\begin{assumption}\label{assumption:boundedness}
$f$ is bounded, i.e. $\norm{f}_\infty < \infty$, and $\inf_{u\in\R}f(u) =: f_{\min}>0. $
\end{assumption}
A simple example of a function satisfying the assumption is the shifted sigmoid. Note that these bounds do not allow taking an exponential function $f$ (or any unbounded function) nor having an absolute refractory period (short interval of time following a spike during which an neuron cannot spike). In other terms, neurons can not be forced to spike in a finite time interval nor forced to stay silent. Nevertheless, since $\norm{f}_\infty$ can be arbitrarily large and $f_{\min}$ can be arbitrarily small, these bounds do not meaningfull alter biological realism.

Finally, to prove that the stationary process exists, we need:
\begin{assumption}\label{assumption:f_deriv}
$f$ is differentiable and $f'$ is bounded. Furthermore, there exists a positive constant $C$ such that $| u f'(u)|\leq C $ for  all $u.$ 
\end{assumption}
This is a purely technical assumption and is rather innocent since $f$ is anyway bounded.

We can now state our main result:
\begin{theorem}\label{theorem:main}
Grant Assumptions~\ref{assumption:constant}--\ref{assumption:f_deriv}. The voltage-structured SPDE~\eqref{eq:SPDE} is stable in variation.
\end{theorem}

The proof is divided into two parts. In the first part, using Meyn-Tweedie theory \cite{MeyTwe93}, we show that the solutions of Eq.~\eqref{eq:SPDE} satisfy a certain recurrence property which then allows us to prove that the associated jump processes couple, using methods from \cite{BreMas96} for nonlinear Hawkes processes. In the second part, we prove the existence of a non-trivial stationary process solving Eq.~\eqref{eq:SPDE}.

In simulations, the simplified model with fixed $\Lambda$, Eq.~\eqref{eq:SPDE}, has a qualitatively similar behavior (from the stability point of view) as the original model of \cite{SchDeg17} where $\Lambda^Z_t$ has an explicit expression in terms of the past $Z$ (see Sec.~\ref{sec:link}). Hence, the proof of Theorem~\ref{theorem:main} provides an important understanding of the role of the `correction term' $\Lambda^Z_t(1-\dots)$ in the original model (Fig.~\ref{fig:1}c--f).

\subsection{Plan of the paper}
First, in Section~\ref{sec:well-posedness}, we prove the well-posedness of the SPDE~\eqref{eq:SPDE} as a measure-valued piecewise deterministic Markov process. The proof of Theorem~\ref{theorem:main} is then presented in Section~\ref{sec:stability}. 

In Section~\ref{sec:link}, we present a concise derivation of the finite-size integral equation~\eqref{eq:stochastic_integral} and a simple simulation algorithm is provided in Section~\ref{sec:simu}. A general simulation algorithm for multiple interacting populations of generalized integrate-and-fire neurons can be found in the Appendix.

\section{Well-posedness}\label{sec:well-posedness}
Although the SPDE~\eqref{eq:SPDE} might look somewhat formal, it can be rigorously formulated in terms of a piecewise deterministic Markov Process (PDMP) taking values in the space $\M_+$ of all positive measures on $\R.$ We endow $\M_+$ with the topology of weak convergence, which makes $\M_+$ Polish. 

Since Assumptions~\ref{assumption:constant} and \ref{assumption:boundedness} are always imposed in the sequel, we will omit their mention. In particular, we will always assume that $f$ is bounded.

For all $\nu \in \M_+$, let us write $(\mathcal{S}(t)\nu)_{t\geq 0}:=(\rho(\cdot,t))_{t\geq 0}$ the solution to the transport equation
\begin{linenomath}\begin{align}
    \partial_t \rho(du,t) - \partial_u\left(\left(\frac{u-\mu}{\taum}\right)\rho(du,t)\right) &= -f(u)\rho(du,t), \qquad \forall (u,t)\in\R\times\R_+^*, \label{eq:transport}\\
    \rho_0 &= \nu.\nonumber
\end{align}\end{linenomath}
\rot{With the notation of \eqref{eq:Phia}, take the flow $ \Phi_{s, t }^0 $ without exterior input, that is, $ a \equiv 0.$ Then we have the explicit representation 
\begin{equation}\label{eq:stnu}
S(t) \nu = \int_\R  \delta_{ \Phi_{0, t }^0 (u) } e^{ - \int_0^t  f ( \Phi^0_{0,r} (u) ) dr } \nu ( du )  . 
\end{equation}
}
$(\mathcal{S}(t))_{t\in\R_+}$ can be seen as a sub-stochastic $\mathcal{C}_0$-semigroup of bounded linear operators on $\M_+$. Moreover, we introduce, for any $ a \in \R_+ $ and any $\nu \in \M_+$, the shifted measure 
\begin{linenomath}\begin{equation*}
\Delta_a \nu : {\mathcal B}(\R) \ni B \mapsto \nu ( (B - a) ). 
\end{equation*}\end{linenomath}
Putting $\rho_0=\nu_0$, we can construct a path-wise solution to Eq.~\eqref{eq:SPDE} following the procedure:
\begin{enumerate}
\item 
We start from an initial value $ \nu_0 \in\M_+$ at time $t=0.$ 
\item
We consider the counting process
\begin{linenomath}\begin{equation*}\label{eq:jump_process}
Z^*_t =  \int_{[0, t ]\times \R_+} \mathbbm{1}_{ z \le N\left[(\mathcal{S}(s)\nu_0)[f] + \Lambda(1-\norm{\mathcal{S}(s)\nu_0})\right]_+} \pi (ds, dz),
\end{equation*}\end{linenomath}
together with its first jump time $ \tau^1 := \inf \{ t \geq 0 : Z^*_t = 1\}$. 
\item
We put $ \rho_t := \mathcal{S}(t)\nu_0$ for all $ t < \tau^1$. 
\item
At time $\tau^1$, we update
\begin{equation}\label{eq:rho_update}
\rho_{\tau^1} := \Delta_{\frac{J}{N}} \left(\mathcal{S}(\tau^1)\nu_0\right) + \frac{1}{N} \delta_0
\end{equation}
and we return to step 1. replacing $\nu_0$ by $\rho_{\tau^1}$ and time $0$ by $\tau^1$.  
\end{enumerate}

\begin{rem}
 This construction provides indeed a PDMP taking values in $ \M_+;$ in between the successive jumps of $Z_t$ only the transport equation acts, and we shall show below that only a finite number of jumps occurs within each finite time interval.  We  have the explicit representation 
\begin{subequations} \label{eq:representation}
\begin{eqnarray}
\rho_t  &= &\int_\R \delta_{\Phi^Z_{0, t }(u) }e^{ - \int_0^t f( \Phi^Z_{0, r}(u) ) dr}  \nu_0 (d u) + \int_{[ 0,  t]} \delta_{ \Phi^Z_{s, t} (0) } e^{ - \int_s^t f( \Phi^Z_{s, r } (0) ) dr } d Z_s, \\
Z_t &=& \frac{1}{N}\int_{[0,t]\times\R_+}\mathbbm{1}_{z\leq N [\rho_{t^-}[f] + \Lambda(1-\norm{\rho_{t^-}})]_+}\pi(ds,dz) .
\end{eqnarray}
\end{subequations}
\rot{In the above formula, the first term on the right hand side of  (\ref{eq:representation}a) corresponds to \eqref{eq:stnu}, except that we have to replace the null exterior input by $ Z$ such that at each jump of $Z,$ the original mass is shifted by $J/ N ,$ according to the jump term $\Delta_{\frac{J}{N}} $ of \eqref{eq:rho_update}. The second term corresponds to the source term $ \frac{1}{N} \delta_0
$ which is added at each jump of $Z$ and then transported by $ S(t) .$ }

The above notion of solution is actually equivalent to the notion of a {\it mild solution} of the SPDE \eqref{eq:SPDE} driven by Poisson noise (see \cite{Fou00} and \cite{Wal86}).
However, since the only underlying noise is Poisson, with finite jump intensity, the notion of a PDMP with values in $\M_+$ seems to be more natural in this context. 
\end{rem}

\begin{rem}
Using the representation Eq.~\eqref{eq:representation}, we can easily make the link between the SPDE~\eqref{eq:SPDE} and the stochastic integral equation~\eqref{eq:stochastic_integral}. Taking the definition of $\bar{A}_t$ in Eq.~\eqref{eq:def_Z}, we have
\begin{linenomath}\begin{align*}
    \bar{A}_t &= [\rho_{t}[f] + \Lambda(1-\norm{\rho_{t}})]_+ \\
    &= \Bigg[\int_{\R_+}f(\Phi_{0,t}^Z(u))e^{ - \int_0^t f( \Phi^Z_{0, r}(u) ) dr} \nu_0(du) + \int_{[ 0,  t]} f(\Phi^Z_{s,t}(0)) e^{ - \int_s^t f( \Phi^Z_{s, r } (0) ) dr } d Z_s\\ &\qquad
    +\Lambda\left( 1 - \int_{\R_+}e^{ - \int_0^t f( \Phi^Z_{0, r}(u) ) dr} \nu_0(du) - \int_{] 0,  t]} e^{ - \int_s^t f( \Phi^Z_{s, r } (0) ) dr } d Z_s\right)\Bigg]_+
\end{align*}\end{linenomath}
(using Eqs.~\eqref{eq:lambda_A}, \eqref{eq:H} and \eqref{eq:H_tilde})
\begin{linenomath}\begin{align*}
    &= \Bigg[H^Z(t) + \int_{[0,t]} \lambda^Z(t|s)S^Z(t|s)dZ_s + \Lambda\left( 1 - \widetilde{H}^Z_t - \int_{[0,t]} S^Z(t|s)dZ_s\right)\Bigg]_+,
\end{align*}\end{linenomath}
showing that Eqs.~\eqref{eq:SPDE} and \eqref{eq:stochastic_integral} are equivalent. Also, since
\begin{linenomath}\begin{equation*}
    \norm{\rho_{t}} = \int_{\R_+}e^{ - \int_0^t f( \Phi^Z_{0, r}(u) ) dr} \nu_0(du) + \int_{] 0,  t]} e^{ - \int_s^t f( \Phi^Z_{s, r } (0) ) dr } d Z_s = \widetilde{H}^Z_t + \int_{[0,t]} S^Z(t|s)dZ_s,
\end{equation*}\end{linenomath}
$\norm{\rho_t}$ is equivalent to the neuronal mass $M_t$ defined in Eq.~\eqref{eq:mass}.
\end{rem}

 \rot{In what follows we study {\it the extended generator} $\mathcal{L}$ of our process, in the sense of Meyn and Tweedie \cite{MeyTwe93}. Extended generators are  defined by the pointwise convergence and the fact that a fundamental martingale property reminiscent of the It\^o formula is verified. For the convenience of the reader we recall its definition: We set
$\mathcal{D}(\mathcal{L} )$ the set of all measurable functions $ g: {\mathcal M}_+ \to \R$ for which there exists a measurable function $h:{\mathcal M}_+\to\R,$ 
 such that $ t \mapsto \mathbb{E}_\nu ( h ( \rho_t) ) $ is continuous in $0,$ and  such that $\forall \nu \in {\mathcal M}_+,$ $\forall t\geq 0,$
\begin{enumerate}
\item $\mathbb{E}_\nu [{g(\rho_t)}] -g(\nu)=\mathbb{E}_\nu{\int_0^t h(\rho_s)ds};$
\item $\mathbb{E}_\nu [\int_0^t|h (\rho_s)|ds] <\infty.$
\end{enumerate}
In this case, we write $ {\mathcal L } g:= h .$} 

On a restricted set of test functions, we can explicitly calculate the extended generator $\mathcal{L}$ of the PDMP described above: For all $\varphi\in\mathcal{C}^1_b(\R)$ (bounded and continuously differentiable functions), for all $\nu\in\M_+$ and using the abuse of notation $\varphi(\nu) := \nu[\varphi]$, we have that
\begin{linenomath}\begin{multline} \label{eq:generator}
    \mathcal{L}\varphi(\nu) = -\int_{\R}\frac{u-\mu}{\taum}\varphi'(u)\nu(du) - \nu[\varphi f] \\
    +N\big[\nu[f] + \Lambda(1-\norm{\nu})\big]_+\left(\int_{\R}\varphi\left(u+\frac{J}{N}\right)\nu(du) + \frac{1}{N}\varphi(0) - \nu[\varphi]\right).
\end{multline}\end{linenomath}

We now show that this process is well-defined. For that sake, let us define, for all $K>0$, the exit time
\begin{equation}\label{eq:exit_time}
T^K:=\inf\{t\geq0\,:\,\norm{\rho_t}>K\} .
\end{equation}

\begin{rem}
The $T^K$ are well-defined stopping times since the sets $\{\nu\in\M_+:\norm{\nu}>K\}$ are the pre-image of $]K,+\infty[$ by the linear form $\mathbf{1}: \M_+ \to \R_+,\,\nu \mapsto \nu[\mathbf{1}]$ and we have endowed $\M_+$ with the topology of weak convergence. For a general treatment of the measurability of hitting times, see \cite{Bas10} and in particular Theorem 2.4 of that article.
\end{rem}

Up to time $ T^K, $ the overall jump intensity of the process is bounded by $ \| f\|_\infty K + \Lambda, $ such that the procedure described above is well-defined up to the explosion  time  of the process
$\zeta := \lim_{K\to+\infty}T^K$.
To show that Eq.~\eqref{eq:SPDE} is well-defined on $\R_+$, we need to prove that the PDMP defined above is \textit{non-explosive} in the sense of \cite{MeyTwe93}, i.e. $\zeta = +\infty$ a. s. We follow the standard `drift condition'-based approach of \cite{MeyTwe93}. Writing $V(\nu) := \norm{\nu} = \nu[\mathbf{1}],\, \forall \nu\in\M_+$, we have

\begin{lemma}[Foster-Lyapunov inequality]\label{lemma:foster_lyapunov}
There exist $K^*>0$, $d >  0$ and $c>0$ such that
\begin{equation} \label{eq:lyapunov}
    \forall \nu\in\M_+, \qquad \mathcal{L}V(\nu) \leq d\mathbbm{1}_{\norm{\nu}\leq K^*} - c(1+V)(\nu).
\end{equation}
\end{lemma}
\begin{proof}
Using Eq.~\eqref{eq:generator} and $V(\nu) = \nu[\mathbf{1}]$, we have $\mathcal{L}V(\nu) = -\nu[f] + \Big[\nu[f] + \Lambda(1-\norm{\nu})\Big]_+$. 

Two cases arise: either $\Big[\nu[f] + \Lambda(1-\norm{\nu})\Big]_+ > 0$, in which case $\mathcal{L}V(\nu)  =   \Lambda\left(1-\norm{\nu}\right) = \Lambda - \Lambda  V ( \nu )$, or $ \Big[\nu[f] + \Lambda(1-\norm{\nu})\Big]_+ = 0$, in which case $\mathcal{L}V(\nu)  =  - \nu[f]  \leq - f_{\min}  V ( \nu )$. 

Whence, $\mathcal{L}V(\nu)\leq \Lambda - (f_{\min} \wedge \Lambda)V(\nu)$. We can adapt the constants to obtain Eq.~\eqref{eq:lyapunov}.
\end{proof}
Arguing as in Theorem 2.1 of \cite{MeyTwe93}, Lemma~\ref{lemma:foster_lyapunov} guarantees that the PDMP is non-explosive. Hence, we have proved the well-posedness of Eq.~\eqref{eq:SPDE}:

\begin{proposition}[Well-posedness]
For all $\nu_0\in\M_+$, there exists a $\M_+$-valued path-wise unique solution to Eq.~\eqref{eq:SPDE} on $\R_+$.
\end{proposition}

\section{Stability}\label{sec:stability}
\subsection{Coupling}
More than non-explosion, the `drift condition'-based method of \cite{MeyTwe93, MeyTwe09} allows us to show that the PDMP~\eqref{eq:SPDE} satisfies a certain `recurrence' property.

For all $K>0$, let us write the hitting time $t_K := \inf\{t\geq 0\,:\,\norm{\rho_t}\leq K\}$ and denote by $\mathbb{E}_{\nu_0}[t_{K}]$ the expected hitting time of the PDMP~\eqref{eq:SPDE} starting in state $\nu_0\in\M_+$ at time $0$.

\begin{lemma}\label{lemma:recurrence}
Take the constant $K^*$ of Lemma~\ref{lemma:foster_lyapunov}. For all $\nu_0\in \M_+$ such that $\norm{\nu_0}>K^*$, $\mathbb{E}_{\nu_0}[t_{K^*}] < +\infty$.
\end{lemma}
\begin{proof}
The proof is standard but we reproduce it here to highlight the fact that it holds even if the space in which the process evolves is not locally compact.

We use $V$ and the constants of Lemma~\ref{lemma:foster_lyapunov}. For any $t>0$ and any $M>K^*$, by Dynkin's formula (see \cite{MeyTwe93}),
\begin{linenomath}\begin{equation*}
\mathbb{E}_{\nu_0}[V(\rho_{t\wedge T^M})] = V(\nu_0) + \mathbb{E}_{\nu_0}\int_0^{t\wedge T^M}\mathcal{L}V(\rho_s)ds \leq V(\nu_0) + dt,
\end{equation*}\end{linenomath}
where $T^M$ is the exit time defined in Eq.~\eqref{eq:exit_time} and where $d$ is given in Eq.~\eqref{eq:lyapunov}. 

Since $ V(\rho_{t\wedge T^M}) \geq M \mathbbm{1}_{ T_M \leq t },$ this implies
$$  \mathbb{P}_{\nu_0} ( T_M \leq t ) \le \frac{V(\nu_0) + dt}{M} .$$ 
Taking $M\to\infty$, by monotone convergence, $  \mathbb{P}_{\nu_0} ( \zeta  \leq t ) = 0, $ which implies non-explosion. 

We now make another use of Dynkin's formula:
\begin{linenomath}\begin{multline*}
    \mathbb{E}_{\nu_0}[V(\rho_{t\wedge t_{K^*} \wedge T^M })] = V(\nu_0) + \mathbb{E}_{\nu_0}\int_0^{t\wedge t_{K^*}\wedge T^M }\mathcal{L}V(\rho_s)ds \\ \leq V(\nu_0) 
    - c\mathbb{E}_{\nu_0}\int_0^{t\wedge t_{K^*}\wedge T^M}(1+V)(\rho_s)ds.
\end{multline*}\end{linenomath}
Whence,
\begin{linenomath}\begin{equation*}
    \mathbb{E}_{\nu_0}\int_0^{t\wedge t_{K^*}\wedge T^M}(1+V)(\rho_s)ds \leq\frac{V(\nu_0) - \mathbb{E}_{\nu_0}[V(\rho_{t\wedge t_{K^*} \wedge T^M})]}{c} \leq \frac{V(\nu_0) - K^*}{c}.
\end{equation*}\end{linenomath}
Taking $t, M\to\infty$, we get, by monotone convergence
\begin{linenomath}\begin{equation*}
\mathbb{E}_{\nu_0}\int_0^{t_{K^*}}(1+V)(\rho_s)ds \leq \frac{V(\nu_0) - K^*}{c}.
\end{equation*}\end{linenomath}
The fact that $\mathbb{E}_{\nu_0}[t_{K^*}]\leq\mathbb{E}_{\nu_0}\int_0^{t_{K^*}}(1+V)(\rho_s)ds$ concludes the proof.
\end{proof}

The definition of stability we use involves the notion of coupling of two processes (see Sec.~\ref{sec:main}). For $\nu_0$ and $\tilde{\nu}_0\in\M_+$, a natural way to couple two processes $\rho$ and $\tilde{\rho}$ following Eq.~\eqref{eq:SPDE} with initial condition $\nu_0$ and $\tilde{\nu}_0$ respectively is to construct them with the same Poisson random measure $\pi$. With this coupling, the associated jump processes $Z$ and $\tilde{Z}_t$ follow, for all $t\geq 0$, 
\begin{linenomath}\begin{equation*}
Z_t := \frac{1}{N} \int_{[0,t]\times\R_+}\mathbbm{1}_{z\leq N[\rho_s[f] + \Lambda(1-\norm{\rho_s})]_+}\pi(ds,dz) , \quad   \tilde{Z}_t := \frac{1}{N}\int_{[0,t]\times\R_+}\mathbbm{1}_{z\leq N[\tilde{\rho}_s[f] + \Lambda(1-\norm{\tilde{\rho}_s})]_+}\pi(ds,dz).
\end{equation*}\end{linenomath}

For all $t\geq 0$, we can now introduce the event
\begin{linenomath}\begin{equation*}
E_t := \{ Z_{t+s}- Z_t =  \tilde Z_{t+s}- \tilde{Z}_t \mbox{ for all } s \geq 0 \}
\end{equation*}\end{linenomath}
on which both jump processes couple after time $t$. With $(\mathcal{F}_t)_{t\geq 0}$ denoting the natural filtration of the coupled process, we have a lower bound on $\mathbb{P}(E_t|\mathcal{F}_t)$:
\begin{lemma} \label{lemma:lower}
For any $K>0$, there exists a constant $ \varepsilon \in ] 0, 1 [ $ such that for all $ t\geq 0$, 
\begin{linenomath}\begin{equation}\label{eq:doeblin}
 {\mathbb{P}} ( E_t | {\cal F}_t ) \geq \varepsilon \mathbbm{1}_{\{\norm{\rho_t} + \norm{\tilde{\rho}_t} \le K \}} .
\end{equation}\end{linenomath} 
\end{lemma}
\begin{proof}
We use the shorthand $\bar{A}[\nu] := [\nu[f] + \Lambda(1-\norm{\nu})]_+$, $\forall \nu \in \M_+$.
Fix any $ t \geq 0 $ such that $\norm{\rho_t} + \norm{\tilde \rho_t} \leq K.$ Write $ \tau^1_t := \inf \{ s > t : (Z_{s}  - Z_t) + ( \tilde{Z}_{s} - \tilde{Z}_t) \geq 1/N \} $ the next jump after time $t$. Noticing that for all $t\leq s<\tau^1_t,$  $ \bar{A}[\rho_{s}] \vee \bar{A}[\tilde{\rho}_{s}] \leq \norm{f}_\infty K + \Lambda$, we clearly have that $t<\tau^1_t$, that is, there is no accumulation of jumps in finite time. 

In what follows, we evaluate the difference $ \bar{A}[\rho_s] - \bar{A}[\tilde{\rho}_s]$, for $t\leq s$.  

We start by considering the difference $\rho_s[f] - \tilde{\rho}_s[f]$, for all $ t \leq s < \tau^1_t$. It is clear that, for all $ t \leq s < \tau^1_t$, 
\begin{linenomath}\begin{equation*}
\rho_s[f] = \int_{\R} \rho_t (du ) f ( \Phi^0_{t,s} (u) ) \exp\left( - \int_t^s f ( \Phi^0_{t,r} (u) ) dr\right) \leq K \norm{f}_\infty e^{-(s-t)f_{\min}},
\end{equation*}\end{linenomath}
where $\Phi^0$ is the flow of the transport equation~\eqref{eq:transport} and where for the inequality, we used the bounds of $f$ given by Assumption~\ref{assumption:boundedness}. Consequently, for all $ t \leq s < \tau^1_t$, $|\rho_s[f] - \tilde{\rho}_s[f]|\leq 2K\norm{f}_\infty e^{-(s-t)f_{\min}}$. Similarly, $|\,\norm{\rho_s}-\norm{\tilde{\rho}_s}|\leq 2K e^{-(s-t)f_{\min}}$.

At the jump time $\tau_t^1,$ two cases arise:
\begin{itemize}
\item \underline{$\tau_t^1$ is an asynchronous jump}, that is, only one of the two processes, say $ Z$, jumps, in which case $\rho$ is shifted to the right by $ J/N, $ and a Dirac mass $ \frac1N \delta_0$ is added (see Eq.~\eqref{eq:rho_update}). Then, for all  $ s \in [\tau_t^1 ,  \tau_t^2[$, where $ \tau_t^2  := \inf \{ s > \tau^1_t : (Z_s - Z_{\tau_t^1}) + (\tilde{Z}_s - \tilde{Z}_{\tau_t^1}) \geq 1/N \}$, we have 
\begin{linenomath}\begin{multline*}
 \rho_s[f] = \int_{\R} \rho_{\tau^1_{t^-} } (du) f( \Phi^0_{\tau_t^1, s} (u + J/ N ) ) \exp\left( - \int_{\tau^1_t}^s f ( \Phi^0_{\tau^1_t,r} (u + J/N) ) dr\right) \\
 + \frac1N f(\Phi^0_{\tau_t^1,s}(0))\exp\left( - \int_{\tau^1_t}^s f ( \Phi^0_{\tau^1_t,r} (0) ) dr\right) ,
\end{multline*}\end{linenomath} 
while 
\begin{linenomath}\begin{equation*}
\tilde{\rho}_s[f] = \int_{\R} \tilde{\rho}_{\tau^1_{t^-} } (du) f( \Phi^0_{\tau_t^1, s} (u) ) \exp\left( - \int_{\tau^1_t}^s f ( \Phi^0_{\tau^1_t,r} (u) ) dr\right).
\end{equation*}\end{linenomath}

As a consequence, 
\begin{linenomath}\begin{multline*}
 | \rho_s[f] - \tilde{\rho}_s[f]| \leq \norm{f}_\infty e^{ - f_{\min} ( s- \tau_t^1)} (\| \rho_{\tau^1_{t^-}}\| + \| \tilde{\rho}_{\tau^1_{t^-}} \| ) + \frac{\norm{f}_\infty}{N} e^{ - f_{\min} ( s - \tau_t^1)}  \\
 \leq 2K\norm{f}_\infty e^{ - f_{\min} ( s- \tau_t^1)} e^{-f_{\min}(\tau^1_t - t)} + \frac{\norm{f}_\infty}{N} e^{ - f_{\min} ( s - \tau_t^1)}    \\
  = 2 K \norm{f}_\infty e^{ - f_{\min} ( s- t )}   + \frac{\norm{f}_\infty}{N} e^{ - f_{\min} ( s - \tau_t^1)}.
\end{multline*}\end{linenomath}

\item \underline{$\tau^1_t$ is a synchronous jump}, in which case we obtain similarly that for all $s \in [\tau_t^1 ,  \tau_t^2[$, 
\begin{linenomath}\begin{equation*}
| \rho_s[f] - \tilde{\rho}_s[f]| \leq 2 K \norm{f}_\infty e^{ - f_{\min} ( s- t )}.
\end{equation*}\end{linenomath}
\end{itemize}

Similar estimates hold for $ |\,\norm{\rho_s} - \norm{\tilde \rho_s}|$. \rot{Since the function $ x \mapsto x_+ $ is Lipschitz with Lipschitz constant $1, $ this implies that 
$$ |\bar{A}[\rho_s] - \bar{A}[\tilde{\rho}_s] |  \le | \rho_s ( f) - \tilde \rho_s ( f) | + \Lambda | \,\norm{\rho_s} - \norm{\tilde \rho_s}|.$$} 

Working iteratively with respect to the successive jump times $ \tau_t^n , n \geq 2, $ and using the above arguments, we deduce that for an appropriate constant $C>0$, for all $t\leq s$, 
\begin{equation}\label{eq:Lip1}
|\bar{A}[\rho_s] - \bar{A}[\tilde{\rho}_s] | \leq C  e^{ - f_{\min} (s-t) } ( \| \rho_t \| + \| \tilde \rho_t \| )  +C  \int_{]t, s ]} e^{ - f_{\min} ( s- r)} d D_r
\end{equation}
where \rot{$(D_s)_{s\geq t}$ is the  process counting the asynchronous jumps of $Z$ and $ \tilde{Z}$. Notice that $(D_s)_{s\geq t}$ has stochastic intensity $ (N  |\bar A[\rho_s] - \bar A[\tilde \rho_s] |  )_{ s \geq t } $.
In particular, the above upper bound implies that on $ [t, \infty [$, $( D_s)_{s\geq t}$ is stochastically upper bounded by a linear Hawkes process, say $ (H_s)_{ s \geq t},$ with self-interaction kernel $h (s) = NCe^{ - \min (f) s} $ and with time inhomogeneous baseline rate $ s\mapsto NC  e^{ - f_{\min} (s- t) } (\| \rho_t \| + \| \tilde \rho_t \| )$. }

\rot{The rest of this proof follows the arguments given in the proof of Theorem~2 of \cite[p.~1581]{BreMas96} together with their Lemma~1. Here are the details of the argument: As a consequence of the above, we obtain the lower bound 
$$  \mathbb{P}(E_t|\mathcal{F}_t) = \mathbb{P}\left( D ( [t, \infty [ ) = 0 \big|\mathcal{F}_t\right)  \geq   \mathbb{P}\left( H ( [t, \infty [ ) = 0 \big|\mathcal{F}_t\right),$$
since $ D$ is stochastically upper bounded by $N.$ But by the structure of the Hawkes process, 
\begin{linenomath}\begin{multline*}
  \mathbb{P}\left( H ( [t, \infty [ ) = 0 \big|\mathcal{F}_t\right)
  = \exp \left(-  \int_t^\infty N C  e^{ - f_{\min} (s- t) } (\| \rho_t \| + \| \tilde \rho_t \|) ds \right)  \\
  =  \exp \left(  - N C  ( \| \rho_t \| + \| \tilde \rho_t \| ) / f_{\min} \right). 
\end{multline*}\end{linenomath}}
 Putting $\varepsilon :=  \exp \left(  - 2 N C K  / f_{\min} \right)$ concludes the proof.
\end{proof}

\begin{theorem}\label{lemma:coupling}
Let $\rho$ and $\tilde{\rho}$ be the coupled processes defined above for initial condition $\nu_0$ and $\tilde{\nu}_0 \in \M_+,$ and write $ \mathbb{E}_{(\nu_0, \tilde{\nu_0})} $ for the associated expectation. Then the associated counting processes $Z$ and $\tilde{Z}$ couple a.s. in finite time, i.e.
\begin{linenomath}\begin{equation*}
    \mathbb{P}\left(\limsup_{t\to+\infty}\left\{(Z_s)_{s\geq t}\neq(\tilde{Z}_s)_{s\geq t}\right\}\right) = 0.
\end{equation*}\end{linenomath}
Moreover, the associated coupling time $ T_c $, defined in \eqref{eq:tc} above, admits exponential moments, that is, there exists a positive constant $ \bar{\lambda} > 0 $ such that for all initial conditions $\nu_0$ and $\tilde{\nu}_0 \in \M_+,$  
\begin{linenomath}\begin{equation}\label{eq:finiteexpmomenttc}
\mathbb{E}_{(\nu_0, \tilde{\nu_0})} [ e^{\bar{\lambda} T_c }] < + \infty .
\end{equation}\end{linenomath}
\end{theorem}
\begin{proof}
The beginning of the proof of this theorem is similar to the Lemma~5 of \cite{BreMas96}. Defining $E_\infty:=\cup_{t=0}^\infty E_t$, $(\mathbb{E}[\mathbbm{1}_{E_\infty}|\mathcal{F}_t])_{t\geq 0}$ is a uniformly integrable martingale and we have $\mathbb{E}[\mathbbm{1}_{E_\infty}|\mathcal{F}_t] \to \mathbbm{1}_{E_\infty}$ a.s.

However, for all $K>0$, we have, by Lemma~\ref{lemma:lower},
\begin{linenomath}\begin{equation*}
    \mathbb{E}[\mathbbm{1}_{E_\infty}|\mathcal{F}_t] = \mathbb{P}(E_\infty|\mathcal{F}_t) \geq \mathbb{P}(E_t|\mathcal{F}_t) \geq \varepsilon \mathbbm{1}_{\{\norm{\rho_t} + \norm{\tilde{\rho}_t} \leq K \}}, \qquad \forall t\geq 0.
\end{equation*}\end{linenomath}
We can easily adapt the proofs of Lemma~\ref{lemma:foster_lyapunov} and Lemma~\ref{lemma:recurrence} to discrete times $n\in \mathbb{N}$ and show that there exists $K^*>0$ such that $\mathbb{P}(\limsup_{n\to\infty}\{\norm{\rho_n} + \norm{\tilde{\rho}_n} \leq K^* \}) = 1$. Hence, $\mathbbm{1}_{E_\infty}\geq \varepsilon$ a.s., which in turn implies that $\mathbb{P}(E_\infty) = 1$. Since the event $E_\infty$ is the complement of the event  $\limsup_{t\to+\infty}\left\{(Z_s)_{s\geq t}\neq(\tilde{Z}_s)_{s\geq t}\right\}$, this concludes the first part of the proof.

The proof of the existence of exponential moments for the coupling time, which is rather classical, is postponed to Appendix~\ref{app:exp}. 
\end{proof}

\subsection{Existence of the stationary process}
We construct a stationary process $Z$ following the lines of \cite{BreMas96}. The main idea is to show that a construction on the whole line $ \R,$ that is, starting from $ t= - \infty $ is feasible. If it is so, then intuitively the constructed process is automatically stationary. More precisely, we have the following theorem. 

\begin{theorem} \label{theorem:existence_stationary_Z}
In addition to the usual assumptions, grant Assumption~\ref{assumption:f_deriv}. Then there exists a unique stationary process $Z$ solving Eq.~\eqref{eq:SPDE}. 
\end{theorem}

\begin{proof}
We only need to show that a stationary process $Z$ exists - uniqueness follows then from the coupling property stated in Lemma \ref{lemma:coupling} above. 

We construct a sequence $ Z^{[n]} $ of jump processes in the following way. For any fixed $n \geq 1, $ let $( \rho^{[n]}, \tilde Z^{[n]})$  be the solution of Eq.~\eqref{eq:SPDE} defined on $ [-n, \infty[, $  starting at time $ - n $ from the initial condition $ \rho^{[n]}_{-n } = \frac1N \delta_0,$ with 
\begin{linenomath}\begin{equation*}
\tilde Z^{[n]}_t =\frac{1}{N}  \int_{[- n ,t]\times\R_+}\mathbbm{1}_{z\leq N\bar{A}^{[n]}_{s^-}}\pi(ds,dz) ,\quad \text{with } \bar{A}^{[n]}_t := \left[\rho^{[n]}_t[f] + \Lambda(1-\| \rho^{[n]}_t \|)\right]_+, \qquad \forall t \geq - n,
\end{equation*}\end{linenomath}
and $\tilde Z^{[n]}_ t \equiv 0 $ for all $ t \le - n .$

In order to obtain a standardized sequence of processes, we put 
\begin{linenomath}\begin{equation*}
  Z^{[n]}_t  :=  \tilde Z^{[n]}_t -  \tilde Z^{[n]}_0. 
\end{equation*}\end{linenomath}
In this way, for all $n,$  $ Z^{[n]} $ is an element of the  Skorokhod space $D(\R,\R)$ with $Z^{[n]}_0=0.$ We shall also consider the associated sequence of processes
\begin{linenomath}\begin{equation*}
 X_s^{[n]} := \rho^{[n]}_s [ f] - \Lambda  \|\rho^{[n]}_s\| , 
\end{equation*}\end{linenomath}
such that the stochastic intensity of $ N Z^{[n]}_s $ is $\lambda^{[n]} (s) :=  N [ X_{s^-}^{ [n]} + \Lambda]_+.$  

{\bf Step 1.} We first show that the family $(Z^{[n]} ,X^{[n]} )_{n\geq 1}$ is tight in the Skorokhod space $D(\R,\R^2).$ 
To do so, we use the criterion of Aldous, see Theorem~{VI.}4.5 of \cite{JacShi13}. It is sufficient to prove that
\begin{itemize}
\item[(a)] for all $ T> 0$, all $\varepsilon >0$,
\begin{linenomath}\begin{equation*}
 \lim_{ \sigma \downarrow 0} \limsup_{n \to \infty } \sup_{ (\tau,\tau') \in P_{\sigma,T}} 
 \mathbb{P} ( |Z_{\tau'}^{[n] } - Z_\tau^{[n]  } | +  |X^{[n]}_{\tau'}  - X^{[n]}_\tau   | > \varepsilon ) = 0,
\end{equation*}\end{linenomath}
where $P_{\sigma,T}$ is the set of all pairs of stopping times $(\tau,\tau')$ such that
$- T \leq \tau \leq \tau'\leq \tau+\sigma\leq T$ a.s.,
\item[(b)] for all $ T> 0$, $\lim_{ K \uparrow \infty } \sup_n 
 \mathbb{P} ( \sup_{ - T \le s \le T } ( |Z_s^{[n] }| + X^{[n]}_s)  \geq K ) = 0$.
\end{itemize}

To check (a), observe that,
\begin{linenomath}\begin{equation*}
  \mathbb{E} [   |Z_{\tau'}^{[n] } - Z_\tau^{[n]  } | ] \le   \frac{1}{N}\mathbb{E} \int_\tau^{\tau + \sigma}\lambda^{[n]} (s) ds \le  \frac{1}{N}\sqrt{2 T \sigma} \sqrt{\sup_{ -T \le s \le T }    \mathbb{E}  \left[ (\lambda^{[n]} (s))^2 \right]} .
\end{equation*}\end{linenomath}
Note that $( \lambda^{[n]} ( s))^2 \le C \| \rho^{[n]}_s\| ^2 + C'$, for some constants $C,C' $ independent of $n$. By similar arguments as in the proof of Lemma~\ref{lemma:foster_lyapunov}, we have that $W (\nu) := \norm{\nu}^2$ satisfies 
\begin{linenomath}\begin{equation} \label{eq:lyapunov2}
    \forall \nu\in\M_+, \qquad \mathcal{L}W(\nu) \leq \alpha  - \beta W(\nu),
\end{equation}\end{linenomath}
for suitable constants $ \alpha, \beta > 0\footnote[2]{See Appendix~\ref{section:appendix_lyapunov2}}.$ Then, it is straightforward to show that \eqref{eq:lyapunov2} implies
\begin{linenomath}\begin{equation*}
  \sup_n \sup_{ -T \le s \le T }  \mathbb{E}[ W ( \rho_s^{[n]} ) ] < \infty , 
\end{equation*}\end{linenomath}
implying (a) for the sequence of processes $ Z^{[n]}.$

We now turn to the study of the sequence of processes $ X^{[n]} .$ We show how to control $\rho^{[n]} [f] ; $ the control of $ \| \rho^{[n]} \| $ is obtained similarly by taking $f\equiv 1$. We fix stopping times $ \tau < \tau' $ and consider the increment $ \rho^{[n]}_{\tau'} [f] - \rho^{[n]}_\tau [f]$ on the event $ Z^{[n]}_{\tau'} -  Z^{[n]}_\tau = 0.$ On this event, 
\begin{linenomath}\begin{equation*}
  \rho^{[n]}_{\tau'}[f] - \rho^{[n]}_\tau[f]  = \int_{\R} \rho^{[n]}_\tau ( du ) \left( f ( \Phi_{\tau, \tau'}^{0}(u)) \exp\left( - \int_\tau^{\tau'} f( \Phi_{\tau,s}^{0} (u)) ds\right)  - f ( u ) \right) . 
\end{equation*}\end{linenomath}
Then, 
\begin{linenomath}\begin{multline*}
 \left| f ( \Phi_{\tau, \tau'}^{0}(u)) \exp\left( - \int_\tau^{\tau'} f( \Phi_{\tau,s}^{0} (u)) ds\right)  - f ( u )\right| \\
 \le | f ( \Phi_{\tau, \tau'}^{0}(u))  - f ( u )|
 + \| f\|_\infty \left| \exp\left( - \int_\tau^{\tau'} f( \Phi_{\tau,s}^{0} (u)) ds\right) - 1 \right| \\
 \le | f ( \Phi_{\tau, \tau'}^{0}(u))  - f ( u )| + \| f\|_\infty (1 - e^{ - \sigma \| f\|_\infty}).
\end{multline*}\end{linenomath} 
Using that $| \Phi_{\tau, \tau'}^{0}(u)  -  u| \le (1 - e^{- \sigma/\taum}) | u - \mu | ,$ Taylor's formula implies 
\begin{linenomath}\begin{equation*}
 | f ( \Phi_{\tau, \tau'}^{0}(u))  - f ( u )| \le |  f' ( \xi )|  (1 - e^{- \sigma/\taum}) | u - \mu |  ,
\end{equation*}\end{linenomath}
where $ \xi \in [ u,  \Phi_{\tau, \tau'}^{0}(u) ] \cup [\Phi_{\tau, \tau'}^{0}(u), u ].$ 

We first produce an upper bound in the case where $ u \geq \mu $ and $ \mu \geq 0 .$  
Since $ |f' (u) |\le C/u $ by Assumption~\ref{assumption:f_deriv} and since $ \xi \geq \Phi_{\tau, \tau'}^{0}(u), $ we have 
\begin{linenomath}\begin{equation}\label{eq:finc}
  | f ( \Phi_{\tau, \tau'}^{0}(u))  - f ( u )| \le C  (1 - e^{- \sigma/\taum}) C_\sigma, 
\end{equation}\end{linenomath}  
where
\begin{linenomath}\begin{equation*}
 C_\sigma := \sup_{ u \geq \mu } \frac{1}{ u e^{ - \sigma/\taum} + \mu (1 - e^{ - \sigma/\taum}) } ( u - \mu ) .
\end{equation*}\end{linenomath}
Moreover, it is clear that, for any $\sigma_0>0$, $ \sup_{ \sigma \le \sigma_0} C_\sigma < \infty$. 

If $ \mu \le 0 $ and $\mu <  u \le 0,$ we use that $ f' ( \xi ) $ is bounded on $ [ \mu , 0 ]$ to obtain Eq.~\eqref{eq:finc}. The case $u<\mu$ is treated analogously.

As a consequence, we get the global upper bound (on the event $Z^{[n]}_{\tau'} -  Z^{[n]}_\tau = 0$):
\begin{linenomath}\begin{equation*}
   \left|\rho^{[n]}_{\tau'} [f]- \rho^{[n]}_\tau [f]\right| \le C  (1 - e^{ - \kappa \sigma })\|  \rho^{[n]}_\tau \| ,\quad\text{with }  \kappa := \|f\|_\infty \vee 1/\taum.
\end{equation*}\end{linenomath}
We conclude the control of $\rho^{[n]}[f]$, on the event $Z^{[n]}_{\tau'} -  Z^{[n]}_\tau =  0$, using the Foster-Lyapunov inequality (Lemma~\ref{lemma:foster_lyapunov}):
\begin{linenomath}\begin{equation*}
  \mathbb{E} \|  \rho^{[n]}_\tau \| \le  \mathbb{E}  \|  \rho^{[n]}_0 \| +  d T, \qquad\text{with $d$ from Eq.~\eqref{eq:lyapunov},}
\end{equation*}\end{linenomath}
and the fact that $\sup_n \mathbb{E}  \|  \rho^{[n]}_0 \|  < \infty$.

To deal with the event $Z^{[n]}_{\tau'} -  Z^{[n]}_\tau > 0,$ observe that 
\begin{linenomath}\begin{equation*}
   \mathbb{E} \left[ \left| \rho^{[n]}_{\tau'} [f]  - \rho^{[n]}_\tau [f] \right| \mathbbm{1}_{\{ Z^{[n]}_{\tau'} -  Z^{[n]}_\tau >  0 \}} \right] \le 
   \|f\|_\infty \mathbb{E} \left[ \left( \| \rho^{[n]}_{\tau'} \| +\| \rho^{[n]}_{\tau} \|  \right)   \mathbbm{1}_{\{ Z^{[n]}_{\tau'} -  Z^{[n]}_\tau >  0 \}} \right].
\end{equation*}\end{linenomath}
Moreover, for any stopping time $\tau$ taking values in between $ - T $ and $T,$ we have
\begin{linenomath}\begin{equation*}
 \mathbb{E}  \left[\| \rho^{[n]}_{\tau} \| \mathbbm{1}_{\{ Z^{[n]}_{\tau'} -  Z^{[n]}_\tau >  0 \}}\right] \le \sqrt{ \mathbb{E}  \| \rho^{[n]}_{\tau} \|^2} \sqrt{\mathbb{P} (Z^{[n]}_{\tau'} -  Z^{[n]}_\tau >  0) }.
\end{equation*}\end{linenomath}
Using similar arguments as above, but now with the Lyapunov function $ W ( \nu ) = \|\nu\|^2 , $ we obtain 
\begin{linenomath}\begin{equation*}
 \sup_n \mathbb{E}  \| \rho^{[n]}_{\tau} \|^2 <  \infty .
\end{equation*}\end{linenomath}
Finally, using the already established control over $Z^{[n]}$, we get that 
\begin{linenomath}\begin{equation*}
 \lim_{ \sigma \downarrow 0 } \sup_n \mathbb{P} (Z^{[n]}_{\tau'} -  Z^{[n]}_\tau >  0) = 0, 
\end{equation*}\end{linenomath}
which concludes the proof of (a).

(b) Let us first observe that $\sup_{ - T \le s \le T }  |Z_s^{[n] }| \le Z_T^{[n]} - Z_{- T }^{[n]}$, and
\begin{linenomath}\begin{equation*}
 \sup_{ - T \le s \le T }  |X_s^{[n] }| \le  C \sup_{ - T \le s \le T }\| \rho_s^{[n] } \| \le C\left(\| \rho_{-T}^{[n] } \|+ Z_T^{[n]} - Z_{-T}^{[n]} \right).
\end{equation*}\end{linenomath}
We can then conclude using the moment estimates established above.

{\bf Step 2.} By tightness we can extract a subsequence $n_k$ such that $ (Z^{[n_k]}, X^{[n_k]} ) $ converges, in $D(\R,\R^2),$ to a limit process that we shall denote $ ( Z, X) .$  We now show that $Z$ is necessarily stationary. For that sake, take a test function $\varphi : D(\R,\R) \to \R_+ $ which is continuous (with respect to the Skorokhod topology), bounded, and which does only depend on $ Z \in D(\R,\R)$ within a finite time interval $[a,b]\subset\R_+$. We have to show that for every $ t \geq 0, $ 
\begin{linenomath}\begin{equation*}
 \mathbb{E} [ \varphi (Z) ] = \mathbb{E} [ \varphi ( \theta_t  Z) ] ,
\end{equation*}\end{linenomath}
where $\theta_t Z $ is the shifted counting process defined by $  (\theta_t Z)_s = Z_{t+s} - Z_t ,$ for all $s \geq 0.$

By weak convergence, we have that 
\begin{linenomath}\begin{equation*}
 \mathbb{E} [ \varphi (Z)]-   \mathbb{E} [\varphi ( \theta_t  Z) ]= \lim_{k \to \infty } \mathbb{E} [ \varphi (Z^{[n_k]})] - \mathbb{E}[ \varphi ( \theta_t  Z^{[n_k]}) ].
\end{equation*}\end{linenomath}

Now we use the coupling property proven in Theorem \ref{lemma:coupling} above. 
For any fixed $k$ and $t$ we realize $Z^{[n_k]} $ and $ \theta_t Z^{[n_k]} $  according to the construction used in the proof of Theorem \ref{lemma:coupling}. 

This means the following. Let $ \pi ( dt, dz) $  be a Poisson random measure on $ \R \times \R_+ $ which has intensity $ dt dz $ on $ \R \times \R_+ .$ We construct $ Z^{[n_k]} $ using the atoms of $ \pi $ within $ [ - n_k, \infty [ \times \R_+ ,$ starting from $ \frac1N \delta_0 $ at time $-n_k.$    Then we choose, independently of $ \pi, $ a random measure $ \tilde \rho_{ - n_k} \sim {\mathcal L} ( \rho^{[n_k]}_{- n_k + t } ) .$ Note that this law does not depend on $n_k;$ it only depends on $t.$  Finally, we realize the process $ \theta_t Z^{[n_k]} $ letting it start at time $-n_k$ from the initial condition  $\tilde \rho_{ - n_k} $ and using the same underlying Poisson random measure $ \pi.$ Let $ T_{coup}^{n_k} $ be the finite coupling time of the two processes. Notice that once again, $ \mathcal L ( T_{coup}^{n_k} ) $ does not depend on $n_k.$ 

Using this coupling, we obtain 
\begin{linenomath}\begin{equation*}
 \left|  \mathbb{E} [ \varphi (Z^{[n_k]})] - \mathbb{E}[ \varphi ( \theta_t  Z^{[n_k]}) ]\right| \le \| \varphi  \|_\infty  \mathbb{P} ( T_{coup}^{n_k} \geq  n_k + a  )   =  \| \varphi  \|_\infty  \mathbb{P} ( T_{coup} > n_k +a  ) \to 0 
\end{equation*}\end{linenomath}
as $ n_k \to \infty , $ implying that $ \mathbb{E} [ \varphi (Z)]-   \mathbb{E} [\varphi ( \theta_t  Z)] = 0.$ Since the test functions $\varphi$ form a separating-class (see Theorem~1.2 in \cite[p.~8]{Bil99}), we have that $Z$ and $\theta_t Z$ have the same law, whence stationarity.

{\bf Step 3.} Now, we verify that the process $Z$, where $Z$ is taken from the stationary limit process $(Z,X)$ constructed above, is a jump process where jumps of size $1/N$ occur with intensity $\lambda_t := N[X_{t^-} + \Lambda]_+$.

To ease the notation, in what follows, we rename the subsequence $n_k$ by $n.$ Using the Skorokhod representation theorem,  we may assume that the above weak convergence is almost sure, for a particular realization of the couples  $ (Z^{[n]},X^{[n]} ) . $ Hence, we know that almost surely, $ (Z^{[n]},X^{[n]} ) \to ( Z, X)$ and $\lambda^{[n]} \to \lambda$. Moreover, let $ \bar Z $ be the process having intensity $ \lambda $ for the same underlying Poisson random measure as (the realization of) $ Z.$ Then, by Fatou's lemma, for any $t  \geq 0, $ 
\begin{linenomath}\begin{equation*}
  \mathbb{E} | Z_t  - \bar Z_t| \le \liminf_n \mathbb{E} | Z^{[n]}_t   - \bar Z_t |  
\le  \frac{1}{N}  \liminf_n \mathbb{E} \int_0^t  | \lambda^{[n]} (s) -  \lambda ( s)  | ds  = 0 ,
\end{equation*}\end{linenomath}
where we used the uniform integrability of the $ \lambda^{[n]} , $ namely that $ \sup_n \sup_{s\in[0,t]} \mathbb{E}[\lambda^{[n]}_s] < \infty.$  The same argument shows that $ \mathbb{E} | Z_t  - \bar Z_t|  = 0 $ for all $ t \le 0.$ Hence $ Z = \bar Z $ almost surely, 
implying that $ Z $ has the limit intensity $ \lambda $. 

{\bf Step 4.} 
Finally, we show that the limit process $Z$ has the right dynamic, i.e. its intensity $\lambda_t$ is equal to $\bar{\lambda}_t$ given by 
\begin{linenomath}\begin{equation}\label{eq:la}
\bar  \lambda_t := N \left[ \sum_{k: T_k < t } \frac1N \exp\left( -\int_{T_k}^{t} f ( \Phi_{T_k,s}^{Z} (0) ) ds \right) ( f (\Phi_{T_k,t}^{Z}(0)) - \Lambda) + \Lambda\right]_+, \qquad \forall t\in\R, 
\end{equation}\end{linenomath} 
where $T_k  $ denote the jump times of $Z$  and $\Phi^{  Z} $ is given in \eqref{eq:Phiz}.

The goal of this step is to show that $ \lambda   \equiv  \bar\lambda .$ Fix some time $t \geq 0 $ and a truncation level $K > 1.$ Since almost surely, $Z$ does not jump at time $t$ nor at time $-K $ for all $ K \geq 1, $ Proposition VI.2.2.1 of \cite{JacShi13} implies that 
 $ Z^{[n]}_t - Z^{[n]}_{-K} \to Z_t - Z_{-K}.$ Therefore, we may choose $ n_K $ be such that  $ Z^{[n]}_t - Z^{[n]}_{-K} = Z_t - Z_{-K}$ for all $ n \geq n_K.$ By the continuity properties of the Skorokhod topology, as $ n \to \infty , $ we have that $ T_k^{[n]} \to T_k $ as $n \to \infty, $ for all $  Z_{-K} \le k \le Z_t$ (Proposition VI.2.2.1 of \cite{JacShi13}). Hence,  
 
\begin{linenomath}\begin{multline*}
 \sum_{k: -K \le T^{[n]}_k < t } \frac1N \exp\left( -\int_{T_k^{[n]}}^{t} f \left( \Phi^{Z^{[n]}}_{T_k^{[n]},s} (0) \right) ds \right) \left( f \left(\Phi_{T_k^{[n]},t}^{Z^{[n]}} (0) \right) - \Lambda\right) \to \\
 \sum_{k: - K \le T_k < t } \frac1N \exp\left( -\int_{T_k}^{t} f ( \Phi^Z_{T_k,s} (0) ) ds \right) ( f (\Phi_{T_k,t}^Z (0) ) - \Lambda) .
\end{multline*}\end{linenomath} 
Notice that the expression on the lhs corresponds to the terms contributing to $X_{t^-}^{[n]}, $ issued by jumps happening after time $ - K.$ Since we know that $ X_t^{[n]} $ converges to $ X_t$ for almost all $t, $ this implies that for all $K,$ 
\begin{linenomath}\begin{multline*}
X_{t^-} =  \sum_{k: - K \le T_k < t } \frac1N \exp\left( -\int_{T_k}^{t} f ( \Phi^Z_{T_k,s} (0) ) ds \right) ( f (\Phi_{T_k,t}^Z (0) ) - \Lambda) \\
 + \lim_{n \to \infty} \sum_{k:  T^{[n]}_k < -K  } \frac1N \exp\left( -\int_{T_k^{[n]}}^{t} f \left( \Phi^{Z^{[n]}}_{T_k^{[n]},s} (0) \right) ds \right) \left( f \left(\Phi_{T_k^{[n]},t}^{Z^{[n]}} (0) \right) - \Lambda\right),
\end{multline*}\end{linenomath}  
where this last limit is necessarily finite. Letting $K \to \infty $ we deduce that 
\begin{linenomath}\begin{multline*}
X_{t^-} =  \sum_{k: T_k < t } \frac1N \exp\left( -\int_{T_k}^{t} f ( \Phi^Z_{T_k,s} (0) ) ds \right) ( f (\Phi_{T_k,t}^Z (0) ) - \Lambda) \\
 + \lim_{K\to\infty}\lim_{n \to \infty} \sum_{k:  T^{[n]}_k < -K  } \frac1N \exp\left( -\int_{T_k^{[n]}}^{t} f \left( \Phi^{Z^{[n]}}_{T_k^{[n]},s} (0) \right) ds \right) \left( f \left(\Phi_{T_k^{[n]},t}^{Z^{[n]}} (0) \right) - \Lambda\right) .
\end{multline*}\end{linenomath}  
Next, we shall prove that 
\begin{linenomath}\begin{equation}\label{eq:TODO}
   \lim_{K\to\infty}\lim_{n \to \infty} \sum_{k:  T^{[n]}_k < -K  } \frac1N \exp\left( -\int_{T_k^{[n]}}^{t} f \left( \Phi^{Z^{[n]}}_{T_k^{[n]},s} (0) \right) ds \right) f \left(\Phi_{T_k^{[n]},t}^{Z^{[n]}} (0) \right) = 0 \quad a.s.,
\end{equation}\end{linenomath}   
a similar argument proving that
\begin{linenomath}\begin{equation*}
   \lim_{K\to\infty}\lim_{n \to \infty} \sum_{k:  T^{[n]}_k < -K  } \frac1N \exp\left( -\int_{T_k^{[n]}}^{t} f \left( \Phi^{Z^{[n]}}_{T_k^{[n]},s} (0) \right) ds \right)\Lambda = 0  \quad a.s.,
\end{equation*}\end{linenomath}  
to obtain that indeed, $\lambda_t = N [ X_{t^-} + \Lambda]_+ = \bar \lambda_t$.

Let us now prove \eqref{eq:TODO}. Using Fatou's lemma, we get 
\begin{linenomath}\begin{multline}\label{eq:TODO3}
  \mathbb{E}    \lim_{K\to\infty}\lim_{n \to \infty} \sum_{k:  T^{[n]}_k < -K  } \frac1N \exp\left( -\int_{T_k^{[n]}}^{t} f \left( \Phi^{Z^{[n]}}_{T_k^{[n]},s} (0) \right) ds \right) f \left(\Phi_{T_k^{[n]},t}^{Z^{[n]}} (0) \right)  \\ 
  \le \liminf_{K\to\infty} \liminf_{n\to\infty}  \mathbb{E}  \sum_{k:  T^{[n]}_k < -K  } \frac1N \exp\left( -\int_{T_k^{[n]}}^{t} f \left( \Phi^{Z^{[n]}}_{T_k^{[n]},s} (0) \right) ds \right) f \left(\Phi_{T_k^{[n]},t}^{Z^{[n]}} (0) \right) .
\end{multline}\end{linenomath}
Using the same arguments as those leading to Eq.~\eqref{eq:Lip1}, we have 
\begin{linenomath}\begin{equation*}
  \sum_{k:  T^{[n]}_k < -K  } \frac1N \exp\left( -\int_{T_k^{[n]}}^{t} f \left( \Phi^{Z^{[n]}}_{T_k^{[n]},s} (0) \right) ds \right) f \left(\Phi_{T_k^{[n]},t}^{Z^{[n]}} (0) \right) \le \|f\|_\infty \| \rho^{[n]}_{- K} \| e^{ - \min ( f) ( t + K ) } . 
\end{equation*}\end{linenomath}
Therefore, the rhs of \eqref{eq:TODO3} is upper bounded by 
\begin{linenomath}\begin{equation*}
 \| f \|_\infty \liminf_{K\to\infty} \liminf_{n\to\infty}  \mathbb{E} ( \| \rho^{[n]}_{- K} \| ) e^{ - \min ( f) ( t + K ) }  = 0, 
\end{equation*}\end{linenomath}
since $\sup_n  \sup_K \mathbb{E} ( \| \rho^{[n]}_{- K} \| ) < \infty$. This concludes the proof.
\end{proof}

\begin{corollary} \label{corollary:existence_stationary}
Under the same assumptions as in Theorem~\ref{theorem:existence_stationary_Z}, there exists a unique stationary process $\{\rho,\widehat{\nu_0}\}$ solving Eq.~\eqref{eq:SPDE}. 
\end{corollary}
\begin{proof}
Taking the process $Z \in D(\R,\R)$ constructed in Theorem~\ref{theorem:existence_stationary_Z} and using the same notations as in Eq.~\eqref{eq:la}, the stationary process $\{\bar{\rho},\bar{\nu}_0\}$ corresponding to $Z$ is simply
\begin{linenomath}\begin{equation*}
    \bar{\nu}_0 =  \sum_{T_k\leq 0} \exp\left(-\int_{T_k}^0 f(\Phi^Z_{T_k,s}(0))ds\right) \frac{1}{N} \delta_{\Phi^Z_{T_k,0}(0)}, 
\end{equation*}\end{linenomath}
and for all $t\geq 0$,
\begin{linenomath}\begin{equation*}
    \bar{\rho}_t = \sum_{T_k\leq t} \exp\left(-\int_{T_k}^t f(\Phi^Z_{T_k,s}(0))ds\right) \frac{1}{N}  \delta_{\Phi^Z_{T_k,t}(0)}.
\end{equation*}\end{linenomath}
\end{proof}

\section{Background on the finite-size population equation}\label{sec:link}

In this section, we first present a concise derivation of the stochastic integral equation~\eqref{eq:stochastic_integral}, which synthesizes the arguments of the original derivation \cite{SchDeg17}. Following the integral equation convention \cite{Ger95,Ger00} and as in \cite{SchDeg17}, we formally put the initial condition at time $-\infty$ and Eq.~\eqref{eq:stochastic_integral} reads as follows: for all $t\in\R$,
\begin{subequations}\label{eq:stochastic_integral_original}
\begin{linenomath}\begin{align}
    dZ_t &= \frac{1}{N}\pi(dt,[0,N\bar{A}_{t^-}]),\\
    \bar{A}_t &= \left[\int_{]-\infty, t]} \lambda^{Z}(t|s) S^{Z}(t|s) dZ_s+ \Lambda_t^Z\left(1-\int_{]-\infty, t]} S^{Z}(t|s)dZ_s\right)\right]_+,
\end{align}\end{linenomath}
where $\pi$ is a  Poisson random measure on $\R\times\R_+$ having Lebesgue intensity and  \rot{$\lambda^Z$ and $S^Z$ are defined by Eq.~\eqref{eq:lambda_A} with replacement Eq.~\eqref{eq:Phiz}. Furthermore,} the time-dependent modulating factor $\Lambda_t^Z$ is given by
\begin{linenomath}\begin{equation}\label{eq:Lambda_original}
    \Lambda^{Z}_t = \frac{\int_{]-\infty, t]}\lambda^{Z}(t|s)\{1-S^{Z}(t|s)\}S^{Z}(t|s)dZ_s}{\int_{]-\infty, t]}\{1-S^{Z}(t|s)\}S^{Z}(t|s)dZ_s}.
\end{equation}\end{linenomath}
\end{subequations}
Note that in the original formulation of the model (see Eqs.~(11) and (12) in \cite{SchDeg17}), the expression for the time-dependent modulating factor $\Lambda^{Z}_t$ involved a `variance function' $v$. Integrating Eq.~(12) in \cite{SchDeg17} gives $v(t|s)=\{1-S^{Z}(t|s)\}S^{Z}(t|s)\dot{Z}_s$. As a consequence, Eq.~(11) in \cite{SchDeg17} can be written as Eq.~\eqref{eq:Lambda_original}, eliminating $v$.

To understand the reasoning behind the derivation of Eq.~\eqref{eq:stochastic_integral_original}, one needs to keep in mind that the goal is to obtain an intensity-based and history-dependent point process (i.e. that only depends on the past $Z$) approximating the empirical population activity of the microscopic model \eqref{eq:network_GL}.

Let $(Z_s)_{s<t}$ denote the past population activity. In terms of $(Z_s)_{s<t}$, the stochastic intensity of the empirical population activity (of the microscopic model), at time $t$, can be expressed as
\begin{linenomath}\begin{equation}\label{eq:stochastic_intensity}
    N \int_{]-\infty,t[}\lambda^Z(t|s)\mathfrak{S}^Z(t|s)dZ_s,
\end{equation}\end{linenomath}
where, for all past spike times $s$, $\mathfrak{S}^Z(t|s)$ denotes  the `microscopic survival processes': if there was a spike at time $s$, $\mathfrak{S}^Z(t|s) = 1$ if the neuron which has fired at time $s$ has not fired a spike in $]s,t[$, and $\mathfrak{S}^Z(t|s) = 0$ if it has. We need to approximate \eqref{eq:stochastic_intensity} by an expression which does not involve the microscopic $\mathfrak{S}^Z(t|s)$, but only the past $Z$. Writing $\Delta\mathfrak{S}^Z(t|s):=  \mathfrak{S}^Z(t|s)-S^Z(t|s)$, we have
\begin{linenomath}\begin{equation} \label{eq:decomposed}
    \int_{]-\infty,t[}\lambda^Z(t|s)\mathfrak{S}^Z(t|s)dZ_s = \int_{]-\infty,t[}\lambda^Z(t|s)S^Z(t|s)dZ_s + \int_{]-\infty,t[}\lambda^Z(t|s)\Delta\mathfrak{S}^Z(t|s)dZ_s.
\end{equation}\end{linenomath}
Note that since the number of neurons $N$ is strictly preserved (in the microscopic model),
\begin{linenomath}\begin{equation}\label{eq:constraint}
    \int_{]-\infty,t[}\Delta\mathfrak{S}^Z(t|s)dZ_s = 1 - \int_{]-\infty,t[}S^Z(t|s)dZ_s.
\end{equation}\end{linenomath}
To replace the microscopic $\Delta\mathfrak{S}^Z(t|s)$ on the RHS of \eqref{eq:decomposed}, we introduce a family of conditionally independent (\rot{conditioned on $Z$}) survival processes $\{(\widehat{\mathfrak{S}}^Z(t'|s))_{t'\geq s}\}_s$ -- one for each past spike time $s<t$ -- defined by
\begin{linenomath}\begin{equation*}
    \widehat{\mathfrak{S}}^Z(t'|s) = \begin{cases} 1 &\text{ if } t'< T_s, \\
    0 &\text{ if } t'\geq T_s
    \end{cases},
\end{equation*}\end{linenomath}
\rot{where $\{T_s\}_{\text{past spike time }s<t}$ are accessory random variables satisfying the following conditions: (i) the variables $\{T_s\}_{\text{past spike time }s<t}$ are conditionally independent given $Z$ and (ii), for all past spike time $s<t$, $T_s$ takes values in $[s, +\infty]$ and satisfies $\mathbb{P}(T_s > t'|Z) = S^Z(t'|s)$, for all $t'\in[s,t[$ ($T_s$ can therefore be interpreted as a `death' time given by the survival $S^Z$). Importantly, the processes $\{(\widehat{\mathfrak{S}}^Z(t'|s))_{t'\in[s,t]}\}_s$ are close but not exactly} equivalent to the microscopic $\{(\mathfrak{S}^Z(t'|s))_{t'\in[s,t]}\}_s$, e.g. the conservation equation~\eqref{eq:constraint} does not hold for the processes $\Delta \widehat{\mathfrak{S}}^Z(t|s):=\widehat{\mathfrak{S}}^Z(t|s)-S^Z(t|s)$. However, the conditional independence of the processes $\{(\widehat{\mathfrak{S}}^Z(t'|s))_{t'\in[s,t]}\}_s$ will allow us to close the system of equations (see below) and this is the reason why they are introduced.

We do the approximation
\begin{linenomath}\begin{equation}
\label{eq:approximation1}
    \int_{]-\infty,t[}\lambda^Z(t|s)\Delta\mathfrak{S}^Z(t|s)dZ_s  \approx \Lambda^Z_t\int_{]-\infty,t[}\Delta\mathfrak{S}^Z(t|s)dZ_s,
\end{equation}\end{linenomath}
where
\rot{\begin{linenomath}\begin{align}
    \Lambda^Z_t &:=  \underset{\Lambda}{\arg \min}\;\mathbb{E}\left[\left(\int_{]-\infty,t[}\left(\lambda^Z(t|s) - \Lambda\right)\Delta\widehat{\mathfrak{S}}^Z(t|s)dZ_s\right)^2 \Bigg|\, Z\right] \nonumber \\
    &=  \underset{\Lambda}{\arg \min}\;\mathbb{E}\left[\int_{]-\infty,t[}\left(\lambda^Z(t|s) - \Lambda\right)^2\Delta\widehat{\mathfrak{S}}^Z(t|s)^2 dZ_s \Bigg|\, Z\right] \nonumber \\
    &=\underset{\Lambda}{\arg \min}\int_{]-\infty,t[}\left(\lambda^Z(t|s) - \Lambda       \right)^2\mathbb{E}\left[\Delta\widehat{\mathfrak{S}}^Z(t|s)^2 \big|\, Z\right]dZ_s. \label{eq:argmin}
\end{align}\end{linenomath}}
Note that in the definition of $\Lambda_t^Z$, Eq.~\eqref{eq:argmin}, we have used $\widehat{\mathfrak{S}}^Z(t|s)$ instead of the microscopic $\mathfrak{S}^Z(t|s)$, which would have defined the minimum conditional mean squared error of the approximation Eq.~\eqref{eq:approximation1}. While this replacement cannot be rigorously justified, it allows us to approximate the conditional mean squared error and, in particular, the position of its minimum. 
Since, \rot{$\mathbb{E}\left[\Delta\widehat{\mathfrak{S}}^Z(t|s)^2 \big|\, Z\right] = \{1 - S^Z(t|s)\}S^Z(t|s)$}, and taking the derivative with respect to $\Lambda$ in \eqref{eq:argmin}, we get
\begin{linenomath}\begin{align*}
    \Lambda^Z_t &=  \frac{\int_{]-\infty,t[}\lambda^Z(t|s)\{1 - S^Z(t|s)\}S^Z(t|s)dZ_s}{\int_{]-\infty,t[}\{1 - S^Z(t|s)\}S^Z(t|s)dZ_s}.
\end{align*}\end{linenomath}
We have obtained an approximation of the stochastic intensity \eqref{eq:stochastic_intensity} which only involves the past $Z$:
\begin{linenomath}\begin{equation*}
    N \int_{]-\infty,t[}\lambda^Z(t|s)\mathfrak{S}^Z(t|s)dZ_s \approx N\left[\int_{]-\infty,t[}\lambda^Z(t|s)S^Z(t|s)dZ_s + \Lambda^Z_t \left(1 - \int_{]-\infty,t[}S^Z(t|s)dZ_s\right)\right]_+.
\end{equation*}\end{linenomath}
(Taking the positive part on the RHS simply guarantees that the intensity is nonnegative.)

In practice, we can deal with the ill-defined initial condition at time $-\infty$ by assuming that $Z_t =0$ for all $t<0$ and $Z_0 = 1$ (all neurons spike at time $0$). Consistently, we also put $\Lambda^Z_0 = 0$. Then, the model Eq.~\eqref{eq:stochastic_integral_original}
can be written

For all $t>0$,
\begin{subequations}\label{eq:infinite_0}
\begin{linenomath}\begin{align}
    Z_t &= 1+\frac{1}{N}\int_{]0,t]\times \mathbb{R}_+}\mathbbm{1}_{z\leq N\bar{A}_{s^-}}\pi(ds,dz), \label{eq:infinite_0_Z}\\
    \bar{A}_t &= \left[\int_{[0, t]} \lambda^{Z}(t|s) S^{Z}(t|s) dZ_s+ \Lambda_t^Z\left(1-\int_{[0,t]} S^{Z}(t|s)dZ_s\right)\right]_+, \label{eq:infinite_0_A}\\
    \Lambda^{Z}_t &= \frac{\int_{[0,t]}\lambda^{Z}(t|s)\{1-S^{Z}(t|s)\}S^{Z}(t|s)dZ_s}{\int_{[0,t]}\{1-S^{Z}(t|s)\}S^{Z}(t|s)dZ_s}, \label{eq:infinite_0_Lambda}
\end{align}\end{linenomath}
\end{subequations}
with the initial condition $Z_0=1$ and $\Lambda^{Z}_0 = 0$. Assuming that the original model Eq.~\eqref{eq:stochastic_integral_original} has the same stability property as the simpler model Eq.~\eqref{eq:SPDE}, this practical choice of initial condition is acceptable as it will be `forgotten' after some time. 

\section{Simulation algorithm}\label{sec:simu}
Here, we present a simple simulation algorithm for Eq.~\eqref{eq:infinite_0}. The algorithm presented below can be easily adapted to the more realistic case of multiple interacting populations for generalized integrate-and-fire neurons \cite{SchDeg17}, as we show in Appendix~\ref{app:multi}.

To ease the notation, here, we drop all the superscripts $Z$. We can rewrite Eq.~\eqref{eq:Phia} and \eqref{eq:lambda_A} as the solution of a SDE: for any $s>0$,
\begin{subequations}\label{eq:ode-single}
\begin{linenomath}\begin{align}
    \label{eq:S-ode-single}
  \frac{dS(t|s)}{dt}&=-\lambda(t|s)S(t|s)\\
  \label{eq:u-t-s-single}
 du(t|s)&=\frac{\mu_t-u(t|s)}{\taum}dt +JdZ_{t}
\end{align}\end{linenomath}
\end{subequations}
with \rot{$\lambda(t|s)=f\bigl(u(t|s)\bigr)$ and} initial conditions $S(s|s)=1$ and $u(s|s)=0$.

\paragraph*{Finite history.}
\label{finite-hist}

For all $t\ge 0$, let us define the free membrane potential $h(t)$ as the solution of
\begin{linenomath}\begin{equation}
  dh_t =\frac{\mu_t-h_t}{\taum}dt +JdZ_{t}
\end{equation}\end{linenomath}
with initial condition $h_0=0$ (cf. Eq.~\eqref{eq:u-t-s-single}). It is clear that for fixed $s>0$, $|u(t|s) - h_t| \to 0$ when $t\to\infty$. In practice, there exists a sufficiently large time $T\gg\taum$ such that for $t-s>T$, the initial condition for Eq.~\eqref{eq:u-t-s-single} will be forgotten and the membrane potential $u(t|s)$ with last reset time $s$ can be well approximated by the free membrane potential $h_t$. We call $T$ the history length. Associated with the free membrane potential is the free hazard rate defined as $\lambda_{free}(t):=f(h_t)$. The free hazard rate can be interpreted as the firing intensity of neurons that have fully recovered from refractoriness because the last spike of those neurons happened before time $t-T$ and thus has been approximately forgotten. For the numerical implementation, it is useful to consider the slightly modified model, in which we use the above approximation, i.e. where $\lambda(t|s)$ is set to $\lambda_{free}(t)$ if $t-s>T$. For the sake of notational simplicity, we will use the same symbols for this approximate model. For $0<t<T$, there is no difference between the approximate and the original model. Hence, the solution of the approximate model is governed by  Eqs.~\eqref{eq:infinite_0} and \eqref{eq:ode-single}. However, for $t>T$, the integrals in Eq.~\eqref{eq:infinite_0_A} and \eqref{eq:infinite_0_Lambda} do not need to be evaluated over the whole history from $0$ to $t$ but reduce to integrals over $]t-T,t]$:
\begin{subequations}\label{eq:finite_0}
\begin{linenomath}\begin{align}
    {\bar{A}}_t &= \left[\int_{]t-T,t]} \lambda(t|s) S(t|s) dZ_s+\lambda_{free}(t) x_t+ \Lambda_t\left(1-\int_{]t-T,t]} S(t|s)dZ_s-x_t\right)\right]_+, \\
    \Lambda_t &= \frac{\int_{]t-T,t]}\lambda(t|s)\{1-S(t|s)\}S(t|s)dZ_s+\lambda_{free}(t)z_t}{\int_{]t-T,t]}\{1-S(t|s)\}S(t|s)dZ_s+z_t}.
\end{align}\end{linenomath}
\end{subequations}
These expressions depend on the additional variables $x_t:=\int_{[0, t-T]} S(t|s) dZ_s$ and $z_t:=\int_{[0, t-T]}\{1-S(t|s)\}S(t|s)dZ_s$ that solve the following SDE's \cite{SchDeg17}:
\begin{subequations}\label{eq:xz_ode}
\begin{linenomath}\begin{align}
    d x_t&=- \lambda_{free}(t)x_t dt+S(t|t-T)d{Z}_{t-T},&x_{T}=0,\\
    d z_t&=-2\lambda_{free}(t) z_t dt+\{1-S(t|t-T)\}S(t|t-T)d{Z}_{t-T},&z_T=0.\label{eq:z-ode}
\end{align}\end{linenomath}
\end{subequations}

\paragraph{Time discretization.}

The model with finite history length Eq.~\eqref{eq:finite_0} with the SDE's~\eqref{eq:ode-single} and \eqref{eq:xz_ode} suggests a straightforward update scheme in discrete time. To this end, we consider an equally-spaced partition of the time-axis with mesh $\Delta t$ and time points $t_{\hat{t}}=\hat{t}\Delta t$, $\hat{t}=0,1,2,\dotsc$. Furthermore, we partition the co-moving history frame $]t-T,t]$ in discrete time points $s_{r,\hat{t}}=(\hat{t}-\mathcal{T}+r)\Delta t$, $r=1,\dotsc,\mathcal{T}$ with $\mathcal{T}=T/\Delta t$. On the discrete time points, we define the following quantities:
\begin{linenomath}\begin{equation*}
  \begin{split}
  n_{r,\hat{t}}:=Z_{s_{r,\hat{t}}+\Delta t}-Z_{s_{r,\hat{t}}},\quad S_{r,\hat{t}}:=S(\hat{t}\Delta t\,|\,s_{r,\hat{t}}),\quad u_{r,\hat{t}}:=u(\hat{t}\Delta t\,|\,s_{r,\hat{t}}),\\
  P_{r,\hat{t}}:=1-\exp\lreckig{-\frac{\Delta t}{2}\lrrund{\lambda(\hat{t}\Delta t\,|\,s_{r,\hat{t}})+\lambda((\hat{t}+1)\Delta t\,|\,s_{r,\hat{t}})}},\\
  h_{\hat{t}}:=h(\hat{t}\Delta t),\quad x_{\hat{t}}:=x(\hat{t}\Delta t),\quad y_{\hat{t}}:=y(\hat{t}\Delta t),\quad z_{\hat{t}}:=z(\hat{t}\Delta t)\\
  \bar P_{\hat{t}}:=1-\exp\lreckig{-\frac{\Delta t}{2}\lrrund{\lambda_{free}(\hat{t}\Delta t)+\lambda_{free}((\hat{t}+1)\Delta t)}}.    
  \end{split}
\end{equation*}\end{linenomath}
Using these quantities, the mesoscopic model can be simulated with the following update rule \cite{SchDeg17}: For $r=1,\dotsc,\mathcal{T}-1$,
\begin{subequations}
  \label{eq:meso-discrte}
\begin{linenomath}\begin{align}
  n_{r,\hat{t}+1}&=n_{r+1,\hat{t}}\\
  S_{r,\hat{t}+1}&=\lrrund{1-P_{r+1,\hat{t}}}S_{r+1,\hat{t}}\\
  u_{r,\hat{t}+1}&=u_{r+1,\hat{t}}+\lrrund{\frac{\mu_{\hat{t}\Delta t}-u_{r+1,\hat{t}}}{\taum}+J\frac{n_{\mathcal{T},\hat{t}}}{\Delta t}}\Delta t\\
  h_{\hat{t}+1}&=h_{\hat{t}}+\lrrund{\frac{\mu_{\hat{t}\Delta t}-h_{\hat{t}}}{\taum}+J\frac{n_{\mathcal{T},\hat{t}}}{\Delta t}}\Delta t\\
  x_{\hat{t}+1}&=\lrrund{1-\bar P_{\hat{t}}}x_{\hat{t}}+S_{1,\hat{t}+1}n_{1,\hat{t}+1}\\
  z_{\hat{t}+1}&=\lrrund{1-\bar P_{\hat{t}}}^2z_{\hat{t}}+P_{\hat{t}}x_{\hat{t}}+\lrrund{1-S_{1,\hat{t}+1}}S_{1,\hat{t}+1}n_{1,\hat{t}+1}
\end{align}\end{linenomath}
with boundary conditions $S_{\mathcal{T},\hat{t}}=1$ and $u_{\mathcal{T},\hat{t}}=0$ for all $\hat{t}>0$, and  
\begin{linenomath}\begin{align}
  n_{\mathcal{T},\hat{t}+1}&=\frac{\xi_{\hat{t}}}{N},\qquad\xi_{\hat{t}}\sim\text{Binomial}(N,\bar n_{\hat{t}}),\\ 
  \bar n_{\hat{t}}&=\bar P_{\hat{t}}x_{\hat{t}}+\sum_{r=2}^{\mathcal{T}}P_{r,\hat{t}}S_{r,\hat{t}}n_{r,\hat{t}}+P_{\Lambda,\hat{t}}\lrrund{1-x_{\hat{t}}-\sum_{r=2}^{\mathcal{T}}S_{r,\hat{t}}n_{r,\hat{t}}},\\
  P_{\Lambda,\hat{t}}&=\frac{\bar P_{\hat{t}}z_{\hat{t}}+\sum_{r=2}^{\mathcal{T}}P_{r,\hat{t}}\lrrund{1-S_{r,\hat{t}}}S_{r,\hat{t}}n_{r,\hat{t}}}{z_{\hat{t}}+\sum_{r=2}^{\mathcal{T}}\lrrund{1-S_{r,\hat{t}}}S_{r,\hat{t}}n_{r,\hat{t}}}.
\end{align}\end{linenomath}
\end{subequations}
The independent, identically distributed binomial random variables $\xi_{\hat{t}}^k$ represent the total number of neurons that fire in the time interval $(\hat{t}\Delta t,(\hat{t}+1)\Delta t]$. Therefore, the empirical population activity, Eq.~\eqref{eq:empir-A}, and the corresponding population rate (intensity) 
are finally obtained as $A_{\hat{t}\Delta t,\Delta t}=n_{\mathcal{T},\hat{t}+1}/\Delta t$ and $\bar A_{\hat{t}\Delta t}=\bar n_{\hat{t}}/\Delta t$, respectively. A pseudo-code implementation of the mesoscopic model, Eq.~\eqref{eq:meso-discrte}, is given in Algorithm~\ref{algo-single}.  \rot{A Julia-code implementation of the extended model (Appendix~\ref{app:multi}, Algorithm~\ref{algo}) is publicly available at the following GitHub link: \url{https://github.com/schwalger/mesodyn-LIF}.}\\

  \begin{algorithm}[H]
    \KwData{External stimulus at grid points $\mu_{\hat{t}\Delta t}$, $\hat{t}=1,\dotsc,t_{sim}$}
    \KwResult{Population activities $A_{\hat{t}\Delta t,\Delta t}$ and rates $\bar{A}_{\hat{t}\Delta t}$, $\hat{t}=1,\dotsc,t_{sim}$}
    \BlankLine

    $\mathcal{T}=\lfloor 5\taum/\Delta t \rfloor +1$\;
    $x=0$, $z=0$, $h=0$\;
    $n_{\mathcal{T}}=1, n_{1:\mathcal{T}-1}=0$\;
    $A_{0,\Delta t} = 1/\Delta t$\;
    $S_{1:\mathcal{T}}=1$, $u_{1:\mathcal{T}}=0$\;
    $\lambda_{free}=f(h)$,  $\lambda_{1:\mathcal{T}}=f(h)$;

    \For{all times $\hat{t}=1,\dotsc,t_{\text{sim}}$}
    {
    $h\leftarrow h+[(\mu_{\hat{t}\Delta t}-h)/\taum+JA_{(\hat{t}-1)\Delta t,\Delta t}]\Delta t$\;
    $P_\lambda=\mathtt{Pfire}(f(h),\lambda_{free})$\;
    $W=P_\lambda x$, $X=x$, $Y=P_\lambda z$, $Z=z$\;
    $x\leftarrow x-W$\;
    $z\leftarrow (1-P_\lambda)^2z+W$\;
    \For{$r=2,\dotsc,\mathcal{T}$}
    {
      $u_{r-1}=u_r+[(\mu_{\hat{t}\Delta t}-u_r)/\taum+JA_{(\hat{t}-1)\Delta t,\Delta t}]\Delta t$\;
      $P_\lambda,\lambda_{r-1}=\mathtt{Pfire}(f(u_{r-1}),\lambda_r)$\;
      $m=S_r n_r$\;
      $v=(1-S_r)m$\;
      $W\leftarrow W+P_\lambda m$\tcp*{$\mathcal{W}:=\int_{[0, t]} \lambda(t|s) S(t|s) dZ_s$}
      $X\leftarrow X+m$\tcp*{$\mathcal{X}:=\int_{[0, t]} S(t|s) dZ_s$}
      $Y\leftarrow Y+P_\lambda v$\tcp*{$\mathcal{Y}:=\int_{[0, t]}\lambda(t|s)\{1-S(t|s)\}S(t|s)dZ_s$}
      $Z\leftarrow Z+v$\tcp*{$\mathcal{Z}:=\int_{[0, t]}\{1-S(t|s)\}S(t|s)dZ_s$}
      $S_{r-1}=(1-P_\lambda)S_r$\;
      $n_{r-1}=n_r$\;
    }
    $x\leftarrow x+S_1n_1$\;
    $z\leftarrow z+(1-S_1)S_1n_1$\;
    if $Z>0$: $P_\Lambda=Y/Z$, else $P_\Lambda=0$\;
    $\bar n=\min(\max(0,W+P_\Lambda(1-X)),1)$\tcp*{expected spike count $N\bar{n}=N\bar A_t\Delta t$}
    draw $n_{\mathcal{T}}=\text{Binomial}(N,\bar n)/N$\;
    $\bar A_{\hat{t}\Delta t}=\bar n/\Delta t$\;
    $A_{\hat{t}\Delta t, \Delta t}=n_{\mathcal{T}}/\Delta t$\;
    }
    \caption{Mesoscopic neuronal population model}\label{algo-single}
  \end{algorithm}
  \begin{function}[H]
    {
      $P_\lambda=(\lambda+\lambda_{old})\Delta t/2$\;
      \textbf{if} $P_\lambda>0.01$ then $P_\lambda\leftarrow 1-e^{-P_\lambda}$\; 
      \textbf{return} $P_\lambda, \lambda$  
    }
    \caption{Pfire($\lambda,\lambda_{old}$)}
  \end{function}

\section{Conclusions}
\label{sec:disc}
We have proven that a simplified version of the model proposed in \cite{SchDeg17} is well-posed and stable in variation in the sense of Brémaud and Massoulié \cite{BreMas96}. The simplified model is a Markov embedding of an intensity-based and history-dependent point process where the history dependence is, loosely speaking, more `nonlinear' than in nonlinear Hawkes processes (in the sense that the past filtering function is updated at each jump event such that even in the argument of the intensity function $f( \cdot)$, the dependence on the past is not linear any more, that is, not given by convolution over the past events). To deal with this difficulty in the proofs, we combined arguments for Markov processes taking values in the space of positive measures and nonlinear Hawkes processes. From this point of view, the finite-size population equation~\eqref{eq:stochastic_integral} is even more `nonlinear', which makes its mathematical analysis challenging. The simplified model and the original model of \cite{SchDeg17} could therefore be seen as examples of general intensity-based and history-dependent point processes, extending nonlinear Hawkes processes. Despite their mathematical complexity, these general point processes are rather practical for applications since they can be efficiently simulated, and, as intensity-based processes, can be easily fitted to empirical data using likelihood-based methods \cite{RenLon20}. We hope that this work will stimulate further mathematical research on these general intensity-based processes, which have already proven to be useful in neuroscience.


\section*{Appendix}
\appendix

\section{Multi-population model}\label{app:multi}

The only difference between the neuron model in Eq.~\eqref{eq:network_GL} and the Generalized integrate-and-fire model considered in \cite{SchDeg17} is the addition of a synaptic filtering kernel $\epsilon$ and an absolute refractory period $\Delta_{\text{abs}}\ge 0$. Accordingly, Eq.~\eqref{eq:network_GL_U} is replaced by 
\begin{linenomath}\begin{equation*}
dU_t^{i,N} = \left[\frac{\mu_t-U_t^{i,N}}{\taum}dt-U_{t^-}^{i,N}dZ_t^{i,N}+\left(\frac{J}{N}\sum_{j=1}^N \int_{]-\infty, t]} \epsilon(t-s) dZ_s^{j,N}\right) dt \right]\mathbbm{1}_{T_t^{i,N}>\Delta_{\text{abs}}},
\end{equation*}\end{linenomath}
where $T_t^{i,N}$ is an additional  ``age''-variable defined by the stochastic dynamics $dT_t^{i,N}=dt-T_{t^-}^{i,N}dZ_t^{i,N}$, which clocks the time elapsed since the last spike of neuron $i$.
Then, the definitions for the hazard rate $\lambda$ and the survival $S$ can be easily adapted replacing $\Phi$ in Eq.~\eqref{eq:Phiz} by
\begin{linenomath}\begin{equation*}
    \Phi^{z}_{s,t}(u):= u e^{-\frac{t-s}{\taum}} + \int_s^te^{-\frac{t-r}{\taum}}\left(\frac{\mu_r}{\taum}+J\int_{]-\infty,r]} \epsilon(r-s')dz_{s'}\right)dr, \qquad \forall u\in\R,
\end{equation*}\end{linenomath}
and replacing $\lambda$ in Eq.~\eqref{eq:lambda_A} by $\lambda^z(t|s)=f(\Phi_{s+\Delta_{\text{abs}},t}^z(0))\mathbbm{1}_{t\ge s+\Delta_{\text{abs}}}$.

As explained in \cite{SchDeg17}, it is straightforward to generalize Eq.~\eqref{eq:infinite_0} (with the aforementioned extensions) to multiple interacting populations. Importantly, the multi-population model allows to coarse-grain microscopic models of large biological networks of neurons, like a cortical column.

Again, we will henceforth drop the superscripts $Z$. Let us consider a system of $K$ interacting (homogeneous) populations, each consisting of $N^1, \dots, N^K$ neurons, with parameters 
\begin{linenomath}\begin{equation*}
    \{N^k, \taum^k, \Delta^k_{\text{abs}}, f^k, \epsilon^k, (\mu_t^k)_{t\geq 0}\}_{k=1\dots,K}
\end{equation*}\end{linenomath}
and average connectivity matrix $\mathbf{J}$, where $J^{kl}$ is the average connection strength from population $l$ to population $k$. The multi-population version of Eq.~\eqref{eq:infinite_0} is 

For all $k = 1,\dots, K$ and $t>0$,
\begin{subequations}\label{eq:multi_infinite}
\begin{linenomath}\begin{align}
  \label{eq:z_t-multi}
    Z_t^k &= 1+\frac{1}{N}\int_{[0,t]\times \R_+}\mathbbm{1}_{z\leq N\bar{A}_{s^-}^k}\pi^k(ds,dz),\\
    \bar{A}_t^k &= \left[\int_{[0,t]} \lambda^k(t|s) S^k(t|s) dZ_s^k+ \Lambda_t^k\left(1-\int_{[0,t]} S^k(t|s)dZ_s^k\right)\right]_+, \label{eq:Abar-multi}\\
    \Lambda_t^k &= \frac{\int_{[0,t]}\lambda^k(t|s)\{1-S^k(t|s)\}S^k(t|s)dZ_s^k}{\int_{[0,t]}\{1-S^k(t|s)\}S^k(t|s)dZ_s^k},  \label{eq:Lam_t-multi}
\end{align}\end{linenomath}
\end{subequations}
with the initial condition $Z_0^1 = \dots = Z_0^K=1$ and $\Lambda_0^1 = \dots = \Lambda_0^K = 0$, where $\{\pi^k\}_{k=1,\dots,K}$ are independent Poisson random measures on $\R_+\times\R_+$ with Lebesgue intensity measure and

\begin{subequations}
\begin{linenomath}\begin{align}
    S^k(t|s) &= \exp\left(-\int_s^t \lambda^k(r|s)dr\right),\label{eq:Sk-int} \\
    \lambda^k(t|s) &= f^k(u^k(t|s))\mathbbm{1}_{t\ge s+\Delta^k},\nonumber\\
    u^k(t|s)&= \mathbbm{1}_{t\ge s+\Delta^k}\int_{s+\Delta^k}^t e^{-\frac{t-r}{\taum^k}}\left(\frac{\mu_r^k}{\taum^k}+\sum_{l=1}^K J^{kl}\int_{[s,r]} \epsilon^k(r-s')dZ_{s'}^l\right)dr.\label{eq:uk-int}
\end{align}\end{linenomath}
\end{subequations}
For simplicity, we have presented here a version of the multi-population model without spike-frequency adaptation nor short-term synaptic plasticity but these features can be included \cite{SchDeg17,SchGer20}.

In the following we choose a delayed expontial synaptic filter $\epsilon^k(t)=\frac{1}{\taus^k}\exp\left(-\frac{t-d^k}{\taus^k}\right)\mathbbm{1}_{t\ge d^k}$, where $\taus^k$ is the synaptic decay time constant and $d^k>0$ denotes the transmission delay associated with the presynaptic population $k$. This choice allows us to rewrite Eq.~\eqref{eq:Sk-int} and \eqref{eq:uk-int} as the solution of a SDE (with delay): for any $s>0$,
\begin{linenomath}\begin{align*}
  \frac{dS^k(t|s)}{dt}&=-\lambda^k(t|s)S^k(t|s),\\
 \taum^k\frac{du^k(t|s)}{dt}&=-u^k(t|s)+\mu_t^k+\taum^k\sum_{l=1}^KJ^{kl}y^l_t,\\
 \taus^kdy^k_t&=-y^k_t dt+dZ_{t-d^k}^k,
\end{align*}\end{linenomath}
with initial conditions $S^k(s|s)=1$, $u^k(s|s)=0$ and $y_0^k=0$.

As in the case for a single population (Section~\ref{sec:simu}), the infinite history of Eq.~\eqref{eq:multi_infinite} can be approximated by a finite history. The method is completely analogous to that described in Section~\ref{sec:simu} except that now, each population $k$ has its own free membrane potential $h^k(t)$ following
\begin{linenomath}\begin{equation*}
  \taum^k\frac{dh^k(t)}{dt}=-h^k(t)+\mu_t^k+\taum^k\sum_{l=1}^KJ^{kl}y^l_t,
\end{equation*}\end{linenomath}
with initial condition $h^k(0)=0$, and its own history length $T^k\gg \taum^k$. 

For the discrete time dynamics, being also completely analogous to the single population case, we get the generalized algorithm:\\

  \begin{algorithm}[H]
    \KwData{External stimulus at grid points $\mu_{\hat{t}\Delta t}^k$, $\hat{t}=1,\dotsc,t_{sim}$, $k=1,\dotsc,K$ }
    \KwResult{Population activities $A_{\hat{t}\Delta t,\Delta t}^k$ and rates $\bar{A}_{\hat{t}\Delta t}^k$, $\hat{t}=1,\dotsc,t_\text{$sim$}$, $k=1,\dotsc,K$}
    \For{all populations $k=1,\dotsc,K$}
    {
      $\mathcal{T}^k=\lfloor(5\taum^k+\Delta^k_{\text{$abs$}})/\Delta t \rfloor+1, \hat{\Delta}^k_{\text{$abs$}} = \lfloor \Delta^k_{\text{$abs$}}/\Delta t\rfloor$, $\hat{d}^k = \lfloor d^k/\Delta t\rfloor$\;
      $x^k=0$, $y^k=0$, $z^k=0$, $h^k=0$\;
      $n_{\mathcal{T}^k}^k=1, n_{1:\mathcal{T}^k-1}^k=0$\;
      $S_{1:\mathcal{T}^k}^k=1$, $u_{1:\mathcal{T}^k}^k=0$\;
      $\lambda_{free}^k=f(h^k)$,  $\lambda_{1:\mathcal{T}^k}^k=f(h^k)$;
    }
    \For{all times $\hat{t}=1,\dotsc,t_{sim}$}
    {
      \textbf{for} \textit{all populations} $k=1,\dotsc,K$ \textbf{do} $I_{syn}^{k}=\sum_{l=1}^KJ^{kl}y^l$\;
      \For{all populations $k=1,\dotsc,K$}
      {
        $h^k\leftarrow h^k+[(\mu_{\hat{t}\Delta t}^k-h^k)/\taum^k+I_{syn}^k]\Delta t$\;
        $P_\lambda,\lambda_{free}^k=\mathtt{Pfire}(f(h^k),\lambda_{free}^k)$\;
        $W=P_\lambda x^k$, $X=x^k$, $Y=P_\lambda z^k$, $Z=z^k$\;
        $x^k\leftarrow x^k-W$\;
        $z^k\leftarrow (1-P_\lambda)^2z^k+W$\;
        \For{$r=2,\dotsc,\mathcal{T}^k-\hat{\Delta}^k_{\text{abs}}$}
        {
          $u_{r-1}^k=u_r^k+[(\mu_{\hat{t}\Delta t}^k-u_r^k)/\taum^k+I_{syn}^k]\Delta t$\;
          $P_\lambda,\lambda_{r-1}^k =\mathtt{Pfire}(f^k(u_{r-1}^k),\lambda^k_r)$\;
          $m=S_r^kn_r^k$\;
          $v=(1-S_r^k)m$\;
          $W\leftarrow W+P_\lambda m$\tcp*{$\mathcal{W}:=\int_{[0, t]} \lambda^k(t|s) S^k(t|s) dZ_s^k$}
          $X\leftarrow X+m$\tcp*{$\mathcal{X}:=\int_{[0, t]} S^k(t|s) dZ_s^k$}
          $Y\leftarrow Y+P_\lambda v$\tcp*{$\mathcal{Y}:=\int_{[0, t]}\lambda^k(t|s)\{1-S^k(t|s)\}S^k(t|s)dZ_s^k$}
          $Z\leftarrow Z+v$\tcp*{$\mathcal{Z}:=\int_{[0, t]}\{1-S^k(t|s)\}S^k(t|s)dZ_s^k$}
          $S_{r-1}^k=(1-P_\lambda)S_r^k$\;
          $n_{r-1}^k=n_r^k$\;
        }
        $x^k\leftarrow x^k+S_1^kn_1^k$\;
        $z^k\leftarrow z^k+(1-S_1^k)S_1^kn_1^k$\;
        \For{time points in refractory period $r=\mathcal{T}^k-\hat{\Delta}^k_{\text{abs}}+1,\dotsc,\mathcal{T}^k$}
        {
          $X\leftarrow X+n_r^k$\;
          $n_{r-1}^k=n_r^k$\;
        }
        if $Z>0$: $P_\Lambda=Y/Z$, else $P_\Lambda=0$\;
        $\bar n=\min(\max(0,W+P_\Lambda(1-X)),1)$\tcp*{expected spike count $N\bar{n}=N\bar A_t^k\Delta t$}
        draw $n_{\mathcal{T}^k}^k=\text{Binomial}(N^k,\bar n)/N^k$\;
        $y^k\leftarrow y^ke^{-\Delta t/\taus^k}+\lrrund{1-e^{-\Delta t/\taus^k}}n_{\mathcal{T}^k-\hat{d}^k}/\Delta t$\;
        $\bar A_t^k=\bar n/\Delta t$\;
        $A_t^k=n_{\mathcal{T}^k}^k/\Delta t$\;
      }
    }
    \caption{Mesoscopic multi-population model with $\Delta_{abs}^k\ge 0$, $d^k\ge 0$, $\taus^k \ge 0$}\label{algo}
  \end{algorithm}

\section{Exponential moments for $T_c$ (end of the proof of Theorem \ref{lemma:coupling})} \label{app:exp}
{\color{black}Introducing $ \bar V ( \nu, \tilde \nu):= \frac12 ( \| \nu \| + \| \tilde \nu \|) $ and $ \bar L $ the generator of the coupled processes $ (\rho_t, \tilde \rho_t),$ we obtain as a direct consequence of \eqref{eq:lyapunov} the control 
$$
\bar L \bar V ( \nu , \tilde \nu ) \le \Lambda - (f_{\min} \wedge \Lambda) \bar V ( \nu , \tilde \nu ) ,    
$$
implying that for any $ 0 < c < f_{\min} \wedge \Lambda ,$ there exists a suitable constant $ K^* $ such that, with $ C := \{ \bar V  \le K^* \} , $
\begin{linenomath}\begin{equation}\label{eq:lyapunovcouple}
\bar L \bar V \le - c \bar V + \Lambda \mathbbm{1}_C .
\end{equation}\end{linenomath}
Fix some $ \delta > 0 $ and introduce the sequence of hitting times 
\begin{linenomath}\begin{equation*}
 T_1 (\delta) = \inf \{ t \geq \delta : (\rho_t, \tilde \rho_t) \in C \} , \quad  T_{n+1} ( \delta ) = \inf \{ t \geq T_n ( \delta ) + \delta : (\rho_t, \tilde \rho_t) \in C \} , \quad n \geq 0 .
\end{equation*}\end{linenomath}
Adapting the arguments of Theorem 3.1 of \cite{DFG} to our frame, we deduce from \eqref{eq:lyapunovcouple} that there exist positive constants $ c_1, \bar{\lambda}$ and $ c ( \delta, \bar{\lambda}) , c_2 ( \delta ) $ with 
\begin{linenomath}\begin{equation*}
 \mathbb{E}_{(\nu_0, \tilde \nu_0)} [ e^{ \bar{\lambda} T_1 ( \delta ) } ] \le c_1 \bar V ( \nu_0, \tilde \nu_0 ) + c_2 ( \delta ) 
\end{equation*}\end{linenomath}
and 
\begin{linenomath}\begin{equation*}
 \mathbb{E}_{(\nu_0, \tilde \nu_0)} [ e^{ \bar{\lambda} (T_{n+1} ( \delta )- T_n ( \delta)  )} ] \le c ( \delta, \bar{\lambda}) \mbox{ for all } n \geq 1.
\end{equation*}\end{linenomath}

Relying on \eqref{eq:doeblin}, we may associate to each  $T_n ( \delta ) $ a Bernoulli random variable $ U_n \sim \mathcal{B} ( \varepsilon), $ independent of $ \mathcal {F}_{T_n( \delta ) } ,$ such that 
\begin{linenomath}\begin{equation*}
 U_n = 1 \mbox{ implies that at time $ T_n (\delta),$ the coupling has succeeded.}
\end{equation*}\end{linenomath}
In particular, 
\begin{linenomath}\begin{equation*}
 T_c \le \inf \{ T_n ( \delta) : U_n = 1 \} 
\end{equation*}\end{linenomath}
and 
\begin{linenomath}\begin{equation*}
  \mathbb{E}_{(\nu_0, \tilde \nu_0)}[ e^{ \bar{\lambda} T_c}] \le \sum_{n=1}^\infty \mathbb{E}_{(\nu_0, \tilde \nu_0)}[ e^{ \bar{\lambda} T_n ( \delta ) } \mathbbm{1}_{ \{ U_1 = \ldots = U_{n-1} = 0\}} ] ,
\end{equation*}\end{linenomath}
for any $ \bar{\lambda} > 0 .$ We are now ready to conclude. Since by monotone convergence, 
\begin{linenomath}\begin{equation*}
 \lim_{\bar{\lambda} \to 0 } \mathbb{E}_{(\nu_0, \tilde \nu_0)} [ e^{ \bar{\lambda} (T_{n+1} ( \delta )- T_n ( \delta)  )} ] = 1 , 
\end{equation*}\end{linenomath}
we choose $ \lambda_c > 0 $ such that for all $ 0 < \bar{\lambda} < \lambda_c,$ 
\begin{linenomath}\begin{equation*}
 \sup_{ n \geq 1 }  \mathbb{E}_{(\nu_0, \tilde \nu_0)} [ e^{2  \bar{\lambda} (T_{n+1} ( \delta )- T_n ( \delta)  )} ] \cdot (1 - \varepsilon ) =: \kappa^2 < 1.
\end{equation*}\end{linenomath}
Using that, by successive conditioning,
\begin{linenomath}\begin{equation*}
 \mathbb{E}_{(\nu_0, \tilde \nu_0)} [ e^{2  \bar{\lambda} T_{n} ( \delta ) } ] \le \mathbb{E}_{(\nu_0, \tilde \nu_0)} [ e^{2  \bar{\lambda} T_{1} ( \delta ) } ]\cdot  \left( \frac{\kappa^2}{1 - \varepsilon}\right)^{n-1}, 
\end{equation*}\end{linenomath}
this implies, using the Cauchy-Schwarz inequality, 
\begin{linenomath}\begin{multline*} \mathbb{E}_{(\nu_0, \tilde \nu_0)}[ e^{ \bar{\lambda} T_c}] \le \sum_{n=1}^\infty \mathbb{E}_{(\nu_0, \tilde \nu_0)}[ e^{ \bar{\lambda} T_n ( \delta ) } \mathbbm{1}_{ \{ U_1 = \ldots = U_{n-1} = 0\}} ] \\
\le \sum_{n=1}^\infty \sqrt{\mathbb{E}_{(\nu_0, \tilde \nu_0)}e^{ 2 \bar{\lambda} T_n ( \delta ) } } (1 - \varepsilon)^{(n-1)/2 } \le \sqrt{\mathbb{E}_{(\nu_0, \tilde \nu_0)}e^{ 2 \bar{\lambda} T_1 ( \delta ) } } \sum_{n=1}^\infty \kappa^{n-1} < \infty , 
\end{multline*}\end{linenomath}
which concludes the proof.
}

\section{Proof of Eq.~\eqref{eq:lyapunov2}}\label{section:appendix_lyapunov2}
Using Eq.~\eqref{eq:generator}, we have $\mathcal{L}W(\nu) = -2 \norm{\nu} \nu[f] + \Big[\nu[f] + \Lambda(1-\norm{\nu})\Big]_+ \left( 2 \norm{\nu}  + \frac{1}{N} \right) . $ Whenever $\Big[\nu[f] + \Lambda(1-\norm{\nu})\Big]_+ > 0, $ this yields, for a suitable constant $C,$  
\begin{linenomath}\begin{equation*}
 \mathcal{L}W(\nu)  \le  -2 W(\nu)+ C ( \; \norm{\nu} +1),
\end{equation*}\end{linenomath}
which implies the claim. The easier case $\Big[\nu[f] + \Lambda(1-\norm{\nu})\Big]_+ = 0 $ follows simply from the fact that $\nu[f]\geq f_{\min}\norm{\nu}$.

\section{Power spectral density}
\label{sec:psd}

In Fig.~\ref{fig:1}b, we have characterized the stationary population activity by the power spectral density (PSD) defined for a wide-sense stationary process $X(t)$ and $f>0$ as \cite{Gar85}
\begin{equation}
\label{eq:psd-definition}
    \tilde{C}_X(f):=\lim_{T\to\infty}\frac{|\tilde{X}_T(f)|^2}{T},\qquad \tilde X_T(f):=\int_{0}^{T}e^{-2\pi i ft}X(t)\,dt.
\end{equation}
For the mesoscopic model, we estimated the PSD from the simulated, empirical population activity $\hat{A}_{t,\mathfrak{h}}^N(t)$, Eq.~\eqref{eq:empir-A} with $\mathfrak{h}=0.001$~s, using the averaged periodogram (Bartlett's method without windowing). Specifically, for the PSD shown in Fig.~\ref{fig:1}, we segmented a $50$~s-long realisation of the empirical population activity (sampled with time step $\mathfrak{h}=0.001$~s) into 50  non-overlapping segments of length $T=1$~s, computed the squared absolute values of the fast Fourier transform for each segment, divided the result by $T$ (as in Eq.~\eqref{eq:psd-definition}) and averaged the resulting periodograms over all 50 segments.

For the microscopic model with $J=0$ (as in Fig.~\ref{fig:1}), the neuronal population consists of $N$ independent renewal processes generated by the LIF model with escape noise. Therefore, the PSD of $A_{t,\mathfrak{h}}^N(t)$ in the limit $\mathfrak{h}\to 0$ is well-known from the renewal formula \cite{Str67I,DegSch14} 
\begin{equation}
    \tilde{C}_{A}(f)=\frac{r}{N}\frac{1-|\tilde{P}_{ISI}(f)|^2}{|1-\tilde{P}_{ISI}(f)|^2}.
\end{equation}
Here, $\tilde{P}_{ISI}(f)=\int_{\mathbb{R}}P_{ISI}(t)e^{-2\pi i ft}\,dt$
is the Fourier transform  of the interspike-interval density of single neurons $P_{ISI}(t)=\lambda^0(t|0)S^0(t|0)\mathbbm{1}_{t\ge 0}$ and $r=\lreckig{\int_0^\infty S^0(t|0)\,dt}^{-1}$ is their firing rate. In Fig.~\ref{fig:1}, these quantities were calculated numerically.

\section*{Author statement}
All authors contributed equally to the present article.

\section*{Acknowledgements}
This work was funded by the Swiss National Science Foundation (grant no.~200020\_184615).  It has also been conducted as part of  the  FAPESP project Research, Innovation and Dissemination Center for Neuromathematics (grant 2013/07699-0) and of the ANR project ANR-19-CE40-0024.

\bibliographystyle{plain}
\bibliography{my}

\begin{thebibliography}{10}

\bibitem{Bas10}
Richard~F. Bass.
\newblock The measurability of hitting times.
\newblock {\em Electron. Commun. Probab.}, 15:99--105, 2010.

\bibitem{BilCai20}
Yazan~N. Billeh, Binghuang Cai, Sergey~L. Gratiy, Kael Dai, Ramakrishnan Iyer,
  Nathan~W. Gouwens, Reza Abbasi-Asl, Xiaoxuan Jia, Joshua~H. Siegle, Shawn~R.
  Olsen, et~al.
\newblock Systematic integration of structural and functional data into
  multi-scale models of mouse primary visual cortex.
\newblock {\em Neuron}, 106(3):388--403, 2020.

\bibitem{Bil99}
Patrick Billingsley.
\newblock {\em Convergence of probability measures}.
\newblock Wiley Series in Probability and Statistics: Probability and
  Statistics. John Wiley \& Sons, Inc., New York, second edition, 1999.

\bibitem{Bre17}
Michael Breakspear.
\newblock Dynamic models of large-scale brain activity.
\newblock {\em Nat. Neurosci.}, 20(3):340--352, 2017.

\bibitem{BreMas96}
Pierre Br\'{e}maud and Laurent Massouli\'{e}.
\newblock Stability of nonlinear {H}awkes processes.
\newblock {\em Ann. Probab.}, 24(3):1563--1588, 1996.

\bibitem{BreMas01}
Pierre Br\'{e}maud and Laurent Massouli\'{e}.
\newblock Hawkes branching point processes without ancestors.
\newblock {\em J. Appl. Probab.}, 38(1):122--135, 2001.

\bibitem{BruHak99}
Nicolas Brunel and Vincent Hakim.
\newblock Fast global oscillations in networks of integrate-and-fire neurons
  with low firing rates.
\newblock {\em Neural Comput.}, 11(7):1621--1671, 1999.

\bibitem{CaiIye16}
Nicholas Cain, Ramakrishnan Iyer, Christof Koch, and Stefan Mihalas.
\newblock The computational properties of a simplified cortical column model.
\newblock {\em PLoS Comput. Biol.}, 12(9):e1005045, 2016.

\bibitem{Che17b}
Julien Chevallier.
\newblock Fluctuations for mean-field interacting age-dependent {H}awkes
  processes.
\newblock {\em Electron. J. Probab.}, 22:49, 2017.

\bibitem{Che17}
Julien Chevallier.
\newblock Mean-field limit of generalized {H}awkes processes.
\newblock {\em Stochastic Process. Appl.}, 127(12):3870--3912, 2017.

\bibitem{CheMel20}
Julien Chevallier, Anna Melnykova, and Irene Tubikanec.
\newblock Theoretical analysis and simulation methods for {H}awkes processes
  and their diffusion approximation.
\newblock {\em arXiv preprint arXiv:2003.10710}, 2020.

\bibitem{ChiGra07}
Anton~V. Chizhov and Lyle~J. Graham.
\newblock Population model of hippocampal pyramidal neurons, linking a
  refractory density approach to conductance-based neurons.
\newblock {\em Phys. Rev. E}, 75(1):011924, 14, 2007.

\bibitem{CorTan20}
Quentin Cormier, Etienne Tanr\'{e}, and Romain Veltz.
\newblock Long time behavior of a mean-field model of interacting neurons.
\newblock {\em Stochastic Process. Appl.}, 130(5):2553--2595, 2020.

\bibitem{CorTan21}
Quentin Cormier, Etienne Tanr{\'e}, and Romain Veltz.
\newblock Hopf bifurcation in a mean-field model of spiking neurons.
\newblock {\em Electronic Journal of Probability}, 26:1--40, 2021.

\bibitem{CosGra20}
Manon Costa, Carl Graham, Laurence Marsalle, and Viet~Chi Tran.
\newblock Renewal in {H}awkes processes with self-excitation and inhibition.
\newblock {\em Adv. in Appl. Probab.}, 52(3):879--915, 2020.

\bibitem{DemGal15}
Anna De~Masi, Antonio Galves, Eva L{\"o}cherbach, and Errico Presutti.
\newblock Hydrodynamic limit for interacting neurons.
\newblock {\em J. Stat. Phys.}, 158(4):866--902, 2015.

\bibitem{DecJir08}
Gustavo Deco, Viktor~K. Jirsa, Peter~A. Robinson, Michael Breakspear, and Karl
  Friston.
\newblock The dynamic brain: from spiking neurons to neural masses and cortical
  fields.
\newblock {\em PLoS Comput. Biol.}, 4(8):e1000092, 2008.

\bibitem{DegSch14}
Moritz Deger, Tilo Schwalger, Richard Naud, and Wulfram Gerstner.
\newblock Fluctuations and information filtering in coupled populations of
  spiking neurons with adaptation.
\newblock {\em Phys. Rev. E}, 90(6):062704, 2014.

\bibitem{DelFou16}
Sylvain Delattre, Nicolas Fournier, and Marc Hoffmann.
\newblock Hawkes processes on large networks.
\newblock {\em Ann. Appl. Probab.}, 26(1):216--261, 2016.

\bibitem{DitLoe17}
Susanne Ditlevsen and Eva L\"{o}cherbach.
\newblock Multi-class oscillating systems of interacting neurons.
\newblock {\em Stochastic Process. Appl.}, 127(6):1840--1869, 2017.

\bibitem{DFG}
Randal Douc, Gersende Fort, and Arnaud Guillin.
\newblock Subgeometric rates of convergence of f-ergodic strong {M}arkov
  processes.
\newblock {\em Stochastic Processes and their Applications}, 119(3):897--923,
  2009.

\bibitem{DumHen16b}
Gr{\'e}gory Dumont, Jacques Henry, and Carmen~O. Tarniceriu.
\newblock Noisy threshold in neuronal models: connections with the noisy leaky
  integrate-and-fire model.
\newblock {\em J. Math. Biol.}, 73(6-7):1413--1436, 2016.

\bibitem{DumPay17}
Gr{\'e}gory Dumont, Alexandre Payeur, and Andr{\'e} Longtin.
\newblock A stochastic-field description of finite-size spiking neural
  networks.
\newblock {\em PLoS Comput. Biol.}, 13(8):e1005691, 2017.

\bibitem{Fou00}
Nicolas Fournier.
\newblock Malliavin calculus for parabolic {SPDE}s with jumps.
\newblock {\em Stochastic Process. Appl.}, 87(1):115--147, 2000.

\bibitem{FouLoe16}
Nicolas Fournier and Eva L\"{o}cherbach.
\newblock On a toy model of interacting neurons.
\newblock {\em Ann. Inst. Henri Poincar\'{e} Probab. Stat.}, 52(4):1844--1876,
  2016.

\bibitem{FriPre19}
Karl Friston, Katrin~H. Preller, Chris Mathys, Hayriye Cagnan, Jakob Heinzle,
  Adeel Razi, and Peter Zeidman.
\newblock Dynamic causal modelling revisited.
\newblock {\em Neuroimage}, 199:730--744, 2019.

\bibitem{GalLoe16}
Antonio Galves and Eva L\"{o}cherbach.
\newblock Modeling networks of spiking neurons as interacting processes with
  memory of variable length.
\newblock {\em J. SFdS}, 157(1):17--32, 2016.

\bibitem{Gar85}
C.~W. Gardiner.
\newblock {\em Handbook of Stochastic Methods}.
\newblock Springer-Verlag, Berlin, 1985.

\bibitem{Ger95}
Wulfram Gerstner.
\newblock Time structure of the activity in neural network models.
\newblock {\em Phys. Rev. E}, 51(1):738, 1995.

\bibitem{Ger00}
Wulfram Gerstner.
\newblock Population dynamics of spiking neurons: fast transients, asynchronous
  states, and locking.
\newblock {\em Neural Comput.}, 12(1):43--89, 2000.

\bibitem{GerKis14}
Wulfram Gerstner, Werner~M Kistler, Richard Naud, and Liam Paninski.
\newblock {\em Neuronal dynamics: From single neurons to networks and models of
  cognition}.
\newblock Cambridge University Press, 2014.

\bibitem{HarShe15}
Kenneth~D. Harris and Gordon~M.G. Shepherd.
\newblock The neocortical circuit: themes and variations.
\newblock {\em Nat. Neurosci.}, 18(2):170--181, 2015.

\bibitem{HeeSta21}
Sophie Heesen and Wilhelm Stannat.
\newblock Fluctuation limits for mean-field interacting nonlinear {H}awkes
  processes.
\newblock {\em Stochastic Process. Appl.}, 139:280--297, 2021.

\bibitem{JacShi13}
Jean Jacod and Albert~N. Shiryaev.
\newblock {\em Limit theorems for stochastic processes}, volume 288 of {\em
  Grundlehren der Mathematischen Wissenschaften [Fundamental Principles of
  Mathematical Sciences]}.
\newblock Springer-Verlag, Berlin, second edition, 2003.

\bibitem{JanRit95}
Ben~H. Jansen and Vincent~G. Rit.
\newblock Electroencephalogram and visual evoked potential generation in a
  mathematical model of coupled cortical columns.
\newblock {\em Biol. Cybern.}, 73(4):357--366, 1995.

\bibitem{LefTom09}
Sandrine Lefort, Christian Tomm, J.-C. Floyd~Sarria, and Carl~C.H. Petersen.
\newblock The excitatory neuronal network of the {C}2 barrel column in mouse
  primary somatosensory cortex.
\newblock {\em Neuron}, 61(2):301--316, 2009.

\bibitem{MeyTwe93}
Sean Meyn and Richard~L. Tweedie.
\newblock Stability of {M}arkovian processes. {III}. {F}oster-{L}yapunov
  criteria for continuous-time processes.
\newblock {\em Adv. in Appl. Probab.}, 25(3):518--548, 1993.

\bibitem{MeyTwe09}
Sean Meyn and Richard~L. Tweedie.
\newblock {\em Markov chains and stochastic stability}.
\newblock Cambridge University Press, Cambridge, second edition, 2009.
\newblock With a prologue by Peter W. Glynn.

\bibitem{PakPer10}
Khashayar Pakdaman, Beno\^{\i}t Perthame, and Delphine Salort.
\newblock Dynamics of a structured neuron population.
\newblock {\em Nonlinearity}, 23(1):55--75, 2010.

\bibitem{PieGal20}
Bastian Pietras, No\'{e} Gallice, and Tilo Schwalger.
\newblock Low-dimensional firing-rate dynamics for populations of renewal-type
  spiking neurons.
\newblock {\em Phys. Rev. E}, 102(2):022407, 23, 2020.

\bibitem{PotDie14}
Tobias~C. Potjans and Markus Diesmann.
\newblock The cell-type specific cortical microcircuit: relating structure and
  activity in a full-scale spiking network model.
\newblock {\em Cereb. Cortex}, 24(3):785--806, 2014.

\bibitem{PozMen15}
Christian Pozzorini, Skander Mensi, Olivier Hagens, Richard Naud, Christof
  Koch, and Wulfram Gerstner.
\newblock Automated high-throughput characterization of single neurons by means
  of simplified spiking models.
\newblock {\em PLoS Comput. Biol.}, 11(6):e1004275, 2015.

\bibitem{RenLon20}
Alexandre Ren{\'e}, Andr{\'e} Longtin, and Jakob~H. Macke.
\newblock Inference of a mesoscopic population model from population spike
  trains.
\newblock {\em Neural Comput.}, 32(8):1448--1498, 2020.

\bibitem{SanKno15}
Paula Sanz-Leon, Stuart~A. Knock, Andreas Spiegler, and Viktor~K. Jirsa.
\newblock Mathematical framework for large-scale brain network modeling in the
  virtual brain.
\newblock {\em Neuroimage}, 111:385--430, 2015.

\bibitem{SchBak18}
Maximilian Schmidt, Rembrandt Bakker, Kelly Shen, Gleb Bezgin, Markus Diesmann,
  and Sacha~J. van Albada.
\newblock A multi-scale layer-resolved spiking network model of resting-state
  dynamics in macaque visual cortical areas.
\newblock {\em PLoS Comput. Biol.}, 14(10):e1006359, 2018.

\bibitem{SchGer20}
Valentin Schmutz, Wulfram Gerstner, and Tilo Schwalger.
\newblock Mesoscopic population equations for spiking neural networks with
  synaptic short-term plasticity.
\newblock {\em J. Math. Neurosci.}, 10:Paper No. 5, 32, 2020.

\bibitem{SchChi19}
Tilo Schwalger and Anton~V. Chizhov.
\newblock Mind the last spike---firing rate models for mesoscopic populations
  of spiking neurons.
\newblock {\em Curr. Opin. Neurobiol.}, 58:155--166, 2019.

\bibitem{SchDeg17}
Tilo Schwalger, Moritz Deger, and Wulfram Gerstner.
\newblock Towards a theory of cortical columns: From spiking neurons to
  interacting neural populations of finite size.
\newblock {\em PLoS Comput. Biol.}, 13(4):e1005507, 2017.

\bibitem{Str67I}
R.~L. Stratonovich.
\newblock {\em Topics in the Theory of Random Noise}, volume~1.
\newblock Gordon and Breach, New York, 1967.

\bibitem{Wal86}
John~B. Walsh.
\newblock An introduction to stochastic partial differential equations.
\newblock In {\em \'{E}cole d'\'{e}t\'{e} de probabilit\'{e}s de
  {S}aint-{F}lour, {XIV}---1984}, volume 1180 of {\em Lecture Notes in Math.},
  pages 265--439. Springer, Berlin, 1986.

\bibitem{WanSch22}
Shuqi Wang, Valentin Schmutz, Guillaume Bellec, and Wulfram Gerstner.
\newblock Mesoscopic modeling of hidden spiking neurons.
\newblock {\em arXiv preprint arXiv:2205.13493}, 2022.

\bibitem{WilCow72}
Hugh~R. Wilson and Jack~D. Cowan.
\newblock Excitatory and inhibitory interactions in localized populations of
  model neurons.
\newblock {\em Biophys. J.}, 12(1):1--24, 1972.

\end{thebibliography}
\end{document}